\documentclass[reqno,a4paper,12pt]{amsart} 

\usepackage{amsmath,amscd,amsfonts,amssymb}
\usepackage{mathrsfs,dsfont}

\usepackage{tikz}
\usetikzlibrary{patterns}
\usepackage{float}
\usepackage{appendix}
\usepackage{caption}
\usepackage{subcaption}
\usepackage{enumerate}
\usepackage{hyperref}

\numberwithin{equation}{section}
\numberwithin{figure}{section}

\def\R{\mathbb{R}}

\def\Z{\mathbb{Z}}

\def\M{\mathcal{M}}

\def\E{\mathcal{E}}

\def\D{\mathcal{D}}
\def\W{\mathcal{W}}

\def\1{\mathds{1}}

\def\dH{\dim_{\mathcal{H}}}

\renewcommand\le{\leqslant}
\renewcommand\ge{\geqslant}
\renewcommand\leq{\leqslant}
\renewcommand\geq{\geqslant}

\renewcommand\hat{\widehat}

\newcommand{\supp}{\operatorname{supp}}

\newcommand{\dist}{\operatorname{dist}}

\theoremstyle{plain}
\newtheorem{thm}{Theorem}[section]
\newtheorem{lem}[thm]{Lemma}

\newtheorem{prop}[thm]{Proposition}

\newtheorem*{claim*}{Claim}
\newtheorem*{thm*}{Theorem}

\theoremstyle{definition}

\newtheorem*{definition*}{Definition}
\newtheorem*{remarks*}{Remarks}
\newtheorem*{remark*}{Remark}

\newenvironment{enumerate-math}
{\begin{enumerate}
\addtolength{\itemsep}{5pt}
}
{\end{enumerate}}

\newenvironment{enumerate-text}
{\begin{enumerate}
\addtolength{\itemsep}{5pt}
}
{\end{enumerate}}

\begin{document}

\title{On the proof of Falconer distance conjecture}

\author[]{Kotaro Inami\textsuperscript{1}}
\author[]{Jianbin Li\textsuperscript{1}}
\author[]{Yudong Liang\textsuperscript{2}}
\author[]{Bochen Liu\textsuperscript{1*}}
\author[]{Francisco Romero Acosta\textsuperscript{1}}

\address[1]{Department of Mathematics, Southern University of Science and Technology,
Shenzhen 518055, China}
\email[Kotaro Inami]{inami@sustech.edu.cn}
\email[Jianbin Li]{leeanalysing@gmail.com}
\email[Bochen Liu]{Bochen.Liu1989@gmail.com}
\email[Francisco Romero Acosta]{romero@sustech.edu.cn}

\address[2]{School of Mathematical Sciences, Beijing Normal University, Beijing 100875, People’s Republic of China}

\email[Yudong Liang]{202321130093@mail.bnu.edu.cn}

\thanks{\textsuperscript{*}Corresponding author.\\
 \textcopyright\ 2026 The authors. All rights reserved.\\
No permission is granted to use this work for training or
fine-tuning artificial intelligence or machine-learning models,
including text and data mining for those purposes, without
prior written consent from the authors.
Rights for such uses are expressly reserved to the extent
permitted by applicable law.
}
\date{}


\begin{abstract}
We are able to catch the main idea from OpenAI's recent release and present a proof of Falconer distance conjecture in our own words. Arguments are judiciously simplified and carefully designed for graduate-level readers. The first draft was hand-typed. LLM then came in for proofreading and making
figures.
\end{abstract}

\maketitle

\setcounter{tocdepth}{1}
\tableofcontents

\section{Introduction}
Falconer distance conjecture \cite{Fal85} is one of the most famous problems in harmonic analysis and geometric measure theory. 
\begin{thm}[Falconer distance conjecture]
	\label{thm}
	Suppose $E\subset\R^d, d\geq 2$, is a compact set of Hausdorff dimension $>d/2$. Then its distance set
	$$\Delta(E):=\{|x-y|: x,y\in E\}$$
	has positive Lebesgue measure. 
\end{thm}
Both assumptions $d\geq 2$ and $\dH E>d/2$ are necessary \cite{Fal85}.

This conjecture was recently claimed to be solved by OpenAI among 721 manuscripts \cite{OAI-Falconer26}. The best currently confirmed result in the plane is due to Guth-Iosevich-Ou-Wang \cite{GIOW18}. An alternative proof was given recently by Gan, Wu, and the corresponding author \cite{GLW26+}, inspired by his earlier work on regular sets \cite{Liu26+}. In higher dimensions, the best currently known results are due to Du-Ou-Ren-Zhang  \cite{DORZ23-1}\cite{DORZ23-2}.

We are not able to verify the proof of OpenAI line by line. However, we are able to catch the main idea from there and present a proof in our own words.

\subsection*{Notation.} $A\lesssim B$ means $A\leq CB$; $A\lessapprox B$ means $A\leq 2^{O(\epsilon)Tk}B$. Denote by $\D_j(\mu)$ the collection of almost disjoint $2^{-Tj}$-cubes $Q$ that intersect $\supp\mu$. Let $\mu_Q:=\mu(Q)^{-1}\mu|_Q$ and $\mu^Q$ be its pushforward under rescaling $x\mapsto 2^{Tj}x$. When $\mu$ is regularized, we use $\mu^{\D_j}$ to denote a representative of $\mu^Q, Q\in\D_j(\mu)$.
\subsection*{Acknowledgements.} This work was partially supported by the National Key R\&D Program of China 2024YFA1015400. The corresponding author would like to thank Professor Yuan Zhou for introducing Yudong Liang into this project.

\section{Preliminaries}\label{sec-prelim}
\subsection{Regularized Frostman measures}
If $E$ is a compact set in the Euclidean space of Hausdorff dimension $\dH E$, then it is well known that for every $s<\dH E$ there exists a finite Borel measure $\mu$ supported on $E$ satisfying
$$\mu(B(x,r))\lesssim r^s,\ \forall\,x\in\R^d, r>0.$$
Such a measure is called a Frostman measure of exponent $s$. We refer to \cite{Mat15} as a reference.

For a general probability measure on the unit cube, the following regularization is due to Keleti-Shmerkin \cite{KS18}.
\begin{prop}
Let $d,T, k\ge1$ be integers, $\epsilon>0$, and $\mathcal D_j$ be the half-open dyadic cubes of side length $2^{-Tj}$ in $[0,1)^d$. Suppose that $\mu$ is a probability measure
on $[0,1)^d$.
Then there exists a finite family $\mathcal X$ of pairwise disjoint sets,
each a union of cubes in $\mathcal D_k(\mu)$, such that:
\begin{enumerate}[(i)]
\item $\displaystyle
 \mu\!\left(\bigcup_{X\in\mathcal X}X\right)\ge1-2^{-\epsilon Tk}.
$
\item $\displaystyle
 \mu(X)\ge2^{-\left(\epsilon+\frac{\log_2(2dT+2)}{T}\right)Tk}, \ \forall\, X\in\mathcal X.
$
\item For every choice $\mathcal{R}(\mu)\in\{\mu_X:X\in\mathcal X \}$,
there exists $\beta=(\beta_1,\ldots,\beta_k)\in[0,d]^k$ such that,
for every $1\le j\le k$ and $Q\in\mathcal D_j$ with
$\mathcal R(\mu)(Q)>0$,
$$
 \mathcal R(\mu)(Q)
 \le2^{-\beta_j T}\mathcal R(\mu)(\widehat Q)
 \le2\mathcal R(\mu)(Q),
$$
where $\widehat Q\in\mathcal D_{j-1}$ is the parent of $Q$.
\end{enumerate}
\end{prop}

In this paper, we call $\mathcal{R}(\mu)$ a $2^{-Tk}$-regularized measure. If, in addition, $\mu$ is Frostman of exponent $s$, then
$$\mathcal{R}(\mu)\big(B(x,r)\big)=\frac{\mu\big(B(x,r)\cap X\big)}{\mu(X)}\lesssim2^{\left(\epsilon+\frac{\log_2(2dT+2)}{T}\right)Tk}r^s\lessapprox r^s,$$
which is still a Frostman measure of exponent $s$ up to a negligible error factor.

Fix any such $\mathcal R(\mu)$ and a corresponding sequence $\beta$.
The branching function $f_{\mathcal R (\mu)}:\{0,1,\dots, k\}\to[0,dTk]$ is defined by
$$ f_{\mathcal R (\mu)}(j):=T\sum_{i=1}^j\beta_i,
 \quad j=1,\ldots,k,\quad f_{\mathcal R (\mu)}(0):=0,$$
which is $dT$-Lipschitz. For every $0\le j\le k$ and $Q\in\mathcal D_j$ with
$\mathcal R(\mu)(Q)>0$, it satisfies
\begin{equation}
	\label{measure-of-cube}
	 2^{-j}2^{-f_{\mathcal R (\mu)}(j)}
 \leq\mathcal R(\mu)(Q)\leq 2^{-f_{\mathcal R (\mu)}(j)}.
\end{equation}
Also denote
$$F_{\mathcal R (\mu)}(j):= f_{\mathcal R (\mu)}(j)-\frac{d}{2}Tj.$$

Notice that, if the frequency is reduced from $2^{Tk}$ down to $2^{T(k-l)}$, then we only need to consider $f_{\mathcal R (\mu)}, F_{\mathcal R (\mu)}$ on $[0,k-l]$, that will be called the branching interval in this paper; if ${\mathcal R (\mu)}$ is reduced to ${\mathcal R (\mu)}^Q, Q\in\D_l(\mu)$, then the associated branching function takes values from $\{f_{\mathcal R (\mu)}(j)-f_{\mathcal R (\mu)}(l), j=l,\dots,k\}$, with $[l,k]$ as the branching interval.

\subsection{A mutli-scale good/bad decomposition}
Let $\mu$ be a Frostman measure on $E$ of exponent $s>d/2$. Then there exist two separated compact subsets $E_1, E_2\subset E$, both of positive measure. To show $\Delta(E)$ has positive Lebesgue measure, it suffices to show
$$\Delta(E_1, E_2):=\{|x_1-x_2|: x_1\in E_1, x_2\in E_2\}$$
has positive Lebesgue measure. Let $\mu_j:=\mu|_{E_j}, j=1,2$. We may assume that every hyperplane has measure zero, otherwise it becomes a problem in lower dimensional spaces. 

The single-scale good/bad decomposition was introduced to study the Lebesgue measure of distance sets by Guth-Iosevich-Ou-Wang (Lemma 3.6 in \cite{GIOW18}), which relies on Orponen's $L^p$-estimates on radial projections \cite{Orp19}. A multi-scale version was developed by the corresponding author (Section 3 in \cite{Liu26+}), who used both Orponen's radial projection estimates and an $L^2$ argument on orthogonal projections. We also refer to \cite{Orp17}\cite{KS18} for earlier work on dimensions of distance sets using the idea of good directions. All these discussions are in the plane, but the idea can be extended to higher dimensions, with Ren's radial projection estimates \cite{Ren23+} in place of Orponen's. In this section, we only state the result for our use, leaving the discussion on its proof to the last section of this paper.

Let $Q_l(x)\in \D_l$ denote the $2^{-Tl}$ dyadic cube containing $x$, and $\omega_{x,y}:=\frac{x-y}{|x-y|}$.

\begin{prop}\label{prop-good-bad}
	Suppose $\mu_1,\mu_2$ are regularized Frostman measures of exponent $s>d/2$ with disjoint support contained in the unit cube, and all hyperplanes have measure zero. Then for all small enough $\epsilon,\eta>0$, there exist $c>0$, a decreasing sequence $\Gamma_k\subset \supp\mu_1\times\supp\mu_2$, subsets $G_k\subset \Gamma_k$, such that
	\begin{itemize}
            \item $\mu_1\times\mu_2(\Gamma_\infty)>0$, where $\Gamma_\infty:=\bigcap_{k=1}^\infty\Gamma_k$.
            \item Each $\Gamma_k$ is a disjoint union of no more than $2^{d\eta Tk}$ Cartesian products $E_1^{k,i}\times E_2^{k,i}$.
	    \item Both  $\mu_1\times\mu_2(\Gamma_k\backslash \Gamma_\infty), \mu_1\times\mu_2(\Gamma_k\backslash G_k)\lesssim 2^{-cTk}$;
		\item for all pairs of points $(x_1,x_2)\in G_{k}$	 and every $0\le l<m\le k$, the following holds.
\begin{enumerate}[1).]
\item If $T$ is a $2^{-Tm}\times\cdots\times 2^{-Tm}\times 2^{-Tl}$ tube parallel to $\omega_{x_1,x_2}$ containing $x_i$, then
\begin{equation}\label{good tube}
 (\mu_i)_{Q_l(x_i)}(T)\lesssim 2^{\epsilon Tk}\cdot 2^{-\frac{d}{2}T(m-l)-\min_{[l,m]} (F_{\mu_i}-F_{\mu_i}(l))}.
\end{equation}
\item If $x\in Q_l(x_i)$ and $S(x)\subset Q_l(x_i)$ denotes the $2^{-Tl}\times\cdots\times 2^{-Tl}\times 2^{-Tm}$ slab  perpendicular to $\omega_{x_1,x_2}$ containing $x$, then
\begin{equation}\label{good slab}
 \int (\mu_i)_{Q_l(x_i)}(S(x))d(\mu_i)_{Q_l(x_i)}(x) \lesssim 2^{\epsilon Tk}\cdot 2^{-T(m-l)-\min_{[l, m]} (F_{\mu_i}-F_{\mu_i}(l))}.
\end{equation}
\end{enumerate}
	\end{itemize}
\end{prop}
In this paper, we call these tubes ``good" and these slabs ``good on average".

\section{The distance measure}\label{sec-distance-measure}
Let $\mu_1, \mu_2$ be Frostman measures on $E_1, E_2$. Define a measure $\nu^{\mu_1, \mu_2}$ on $\Delta(E_1, E_2)$ by
$$\int f(t)\,d\nu^{\mu_1, \mu_2}(t):=\iint_{(x_1, x_2)\in\Gamma_\infty} f(|x_1-x_2|)\,|x_1-x_2|^{M-\frac{d-1}{2}}d\mu_1(x_1)d\mu_2(x_2),$$
where $\Gamma_\infty$ is from Proposition \ref{prop-good-bad} and $M=M_{d,s}\in\Z_+$ will be clarified in Section \ref{sec-Mattila-integral}.

Our plan is to show the Radon-Nikodym derivative $d\nu/dt\in L_{loc}^1(0,\infty)$, then $\supp\nu$ must have positive Lebesgue measure. More precisely, let $\widehat{P_k\nu}(\xi)=\hat{\nu}(\xi)\phi(\xi/2^{Tk}), k\geq 1$, be the Littlewood-Paley piece, with $P_0\nu:=\nu-\sum_{k\geq 1}P_k \nu$. If there exists $\rho>0$ such that for each regularized piece $\nu^{\mathcal{R}_k(\mu_1), \mathcal{R}_k(\mu_2)}$, one can find a measure $\nu_k^{\mathcal{R}_k(\mu_1), \mathcal{R}_k(\mu_2)}$ satisfying the total variation estimate
\begin{equation}\label{L-1-estimate}\|\nu^{\mathcal{R}_k(\mu_1), \mathcal{R}_k(\mu_2)}-\nu_k^{\mathcal{R}_k(\mu_1), \mathcal{R}_k(\mu_2)}\|\leq \mathcal{R}_k(\mu_1)\times \mathcal{R}_k(\mu_2)(\Gamma_k\backslash (G_k\cap\Gamma_{\infty}))\end{equation}
and the $L^2$ estimate
\begin{equation}\label{goal-L-2}\|P_k\nu_k^{\mathcal{R}_k(\mu_1), \mathcal{R}_k(\mu_2)}\|_{L^2}\lesssim 2^{-\rho Tk},\end{equation}
then
$$\nu^{\mu_1, \mu_2}=P_0\nu^{\mu_1, \mu_2}+\sum_{k\geq 1}P_k(\nu^{\mu_1, \mu_2}-\nu_k^{\mu_1, \mu_2})+\sum_{k\geq 1}P_k\nu_k^{\mu_1, \mu_2} $$
is convergent in $L_{loc}^1$, as desired.

From now we only consider regularized Frostman measures. So we write $\mu_j$ for $\mathcal{R}_k(\mu_j), j=1,2$, and $\nu_k$ for $\nu_k^{\mathcal{R}_k(\mu_1), \mathcal{R}_k(\mu_2)}$, that will not lead to any confusion.

Now we define the auxiliary measure $\nu_k$. Recall $G_k$ from Proposition \ref{prop-good-bad}.

Let $0<\epsilon\ll\delta\ll 1$ be constants depending on $s,d$ that will be clarified later.

For each $1\leq j\leq k$ and each $Q\in\D_j(\mu_1)$, we say a direction $\omega\in S^{d-1}$ is good for $Q$ if there exist $x_1\in Q$ and $x_2\in\supp\mu_2$ such that
$$\omega_{x_1, x_2}=\pm\omega \quad \text{and} \quad (x_1,x_2)\in G_k.$$
It is similar for ``$\omega$ is good for $Q\in\D_j(\mu_2)$". Then, for $m<j$, take $\Theta_{Q,m}\in C_0^\infty(S^{d-1})$, $\Theta_{Q,m}\leq 1$, to be supported on the $2\cdot 2^{-(1-\epsilon)T(j-m)}$-neighborhood of
\begin{equation}\label{def-Theta}\{\omega\in S^{d-1}: \omega \text{ is good for } Q\},\end{equation}
 and equal to $1$ on the $2^{-(1-\epsilon)T(j-m)}$-neighborhood of \eqref{def-Theta}. Moreover, we can take $\Theta_{Q,m}$ to satisfy
\begin{equation}\label{derivative-Theta}\|\partial_{\omega}^n\Theta_{Q,m}\|_{L^\infty}\lesssim_n 2^{(1-\epsilon)nT(j-m)}.\end{equation}

When $\mu_1, \mu_2$ are treated separately, we just write $\mu$ to simplify the notation.

In this paper, we only consider branching intervals of length $>10\delta k$.

For each $[h_1,h_2]\subset[0,k]$ of length $>10\delta k$, let
\begin{equation}\label{def-j-h-1-h-2}j_{[h_1, h_2]}:=\min \{j\in [h_1+\delta k, h_2-\delta k]: F_\mu(j)=\min_{[h_1+\delta k, h_2-\delta k]} F_\mu\}\end{equation}
and denote the distance from $j_{[h_1, h_2]}$ to the boundary of $[h_1, h_2]$ by 
\begin{equation}\label{def-d-h}d_{[h_1, h_2]}:=\min\{j_{[h_1,h_2]}-h_1, h_2-j_{[h_1, h_2]}\}\in[\delta k, (h_2-h_1)/2].\end{equation}
For $0\leq l_1\leq h<l_2\leq k$, we will also need
\begin{equation}\label{def-d-max}d^{[l_1,l_2]}_h:=\max\{d_{[h_1,h_2]}: [h_1, h_2]\subset [l_1,l_2], j_{[h_1,h_2]}>h\}\end{equation}
and
\begin{equation}\label{def-j-max}j^{[l_1,l_2]}_h:=\max\{j_{[h_1,h_2]}>h: [h_1, h_2]\subset [l_1,l_2], d_{[h_1,h_2]}=d^{[l_1,l_2]}_h\}.\end{equation}

Then we define microlocal cutoff functions
\begin{equation}\label{def-b-k}b^{[l_1,l_2]}_h(x,\omega):=\prod_{\substack{[h_1,h_2]\subset [l_1,l_2]\\0<j_{[h_1,h_2]}-m\leq d_{[h_1, h_2]}\\j_{[h_1, h_2]}>h}}\left(\sum_{Q\in\mathcal{D}_{j_{[h_1, h_2]}}(\mu)} \chi_Q(2^{-Th}x)\Theta_{Q,m}(\omega)\right),\end{equation}
and
\begin{equation}\label{def-b-k-initial}b_k(x,\omega):= \left(\sum_{Q\in\mathcal{D}_{k_0}(\mu)} \chi_Q(x)\Theta_{Q,0}(\omega)\right)\cdot b^{[0,k]}_0(x,\omega),\end{equation}
where $k_0\leq k/3$ only depends on $F_\mu$ that will be chosen at the beginning of Section \ref{sec-proof-main}. The factor for $k_0$ will be removed immediately after the very first step of our reduction. Notice that
\begin{equation}\label{derivative-b}\|\partial_{\omega}^nb^{[l_1,l_2]}_h\|_{L^\infty}\lesssim_n2^{(1-\epsilon)nT\cdot d^{[l_1,l_2]}_h} \quad \text{and}\quad \|\partial_{\omega}^nb_k\|_{L^\infty}\lesssim_n2^{(1-\epsilon)nTk/2}.\end{equation}
Two observations on the construction of $b^{[l_1,l_2]}_h$ will help below when the scale changes:
\begin{enumerate}[(i)]
    \item $b^{[l_1,l_2]}_h$ decreases when $[l_1,l_2]$ is enlarged and
    $$b^{[l_1,l_2]}_h(2^{-T(l-h)}x,\omega)\leq b^{[l_1,l_2]}_l(x,\omega), \ \forall\,h\leq l<l_2.$$
    \item if we only consider $x\in Q_0\in\D_l$, $l\geq h$, then factors for $j_{[h_1, h_2]}\leq l$ in \eqref{def-b-k} are functions on $\omega$ independent in $x$.
\end{enumerate}

Finally, define
\begin{equation}\label{def-nu-k}\begin{aligned}&\int f(t)\,d\nu_k(t)\\:=&\iint_{(x_1, x_2)\in\Gamma_k} f(|x_1-x_2|)|x_1-x_2|^{M-\frac{d-1}{2}}b_k(x_1,\omega_{x_1,x_2})b_k(x_2,\omega_{x_1,x_2})d\mu_1(x_1)d\mu_2(x_2)\\
=& \sum_{i}\int_{E_2^{k,i}}\int_{E_1^{k,i}} f(|x_1-x_2|)|x_1-x_2|^{M-\frac{d-1}{2}}b_k(x_1,\omega_{x_1,x_2})b_k(x_2,\omega_{x_1,x_2})d\mu_1(x_1)d\mu_2(x_2)
 \end{aligned}\end{equation}
Since one of $b_{k}(x_1,\omega_{x_1,x_2}), b_{k}(x_2,\omega_{x_1,x_2})\neq 1$ means $(x_1, x_2)\notin G_k$, we have \eqref{L-1-estimate} due to Proposition \ref{prop-good-bad}.

The rest of this paper is devoted to proving \eqref{goal-L-2}. The proof strategy is to shorten the length of the branching interval from $k$ to $10\delta k$. Lemma \ref{lem-left-to-right} cuts the branching interval from the left, which also provides an absolute power-saving as an initial step. Lemma \ref{lem-proj} cuts the branching interval from the right. However, these two lemmas work on different quantities, $\M$ defined in \eqref{def-M} and $\E$ defined in \eqref{def-E}, so we also need Lemma \ref{lem-high-low} and Lemma \ref{lem-from-M-to-E}  for transitions.

\section{The microlocal Mattila integral}\label{sec-Mattila-integral}
As there are $\lesssim 2^{d\eta Tk}$ terms in the sum in \eqref{def-nu-k} and $\eta$ can be arbitrarily small, we only consider one of them and omit the integral domain $E_1^{k,i}\times E_2^{k,i}$. 
\begin{lem}[The microlocal Mattila integral]\label{lem-Mattila}
\begin{equation}\label{def-M}\begin{aligned}\left(\int_{r\approx 2^{Tk}}|\widehat{\nu_k}(r)|^2dr\right)^2\lessapprox&\prod_{j=1,2}\int_{r\approx 2^{Tk}}\left(\int \left|\int e^{-2\pi irx\cdot\omega}b_k(x,\omega)\,d\mu_j(x)\right|^2d\omega\right)^2r^{d-1}dr\\:=&\prod_{j=1,2}\mathcal{M}_k(\mu_j).\end{aligned}\end{equation}
\end{lem}
A remark is, the estimate in Lemma \ref{lem-Mattila} still holds if the integral domain over $r$ is slightly extended to $2^{(1-\epsilon^2)Tk}\lesssim r\lesssim 2^{Tk}$. This extra $\epsilon^2$-power in the defition of $\M_k$ will not influence our computation below, so we omit it sometimes. Also it would be equivalent if we add a smooth cutoff $\phi(r/2^{Tk})$ into the integral, where $\phi\in C_0^\infty([2^{-T},2^T])$.

We call $\mathcal{M}_k(\mu)$ the microlocal Mattila integral, because it becomes the classical Mattila integral \cite{Mat87} without the microlocal cutoff:
$$\int_{r\approx 2^{Tk}}\left(\int_{S^{d-1}} |\hat{\mu}(r\omega)|^2 d\omega\right)^2r^{d-1}dr.$$

To prove Lemma \ref{lem-Mattila}, we need a ``reversed" version of the stationary phase expansion that has its own interest.
\begin{lem}[Reversed stationary phase expansion]\label{lem-reverse-stationary-phase}
	 Suppose $a\in C_0^\infty(S^{d-1})$ is supported on a local chart and $\omega_{x_1,x_2}\in \supp a$. Then for $M\in\Z_+$,
\begin{multline*}
	e^{-2\pi i r|x_1-x_2|}|x_1-x_2|^{M-\frac{d-1}{2}}a(\omega_{x_1,x_2})\\= C_d\,\sum_{j=0}^M r^{\frac{d-1}{2}}\int e^{-2\pi ir\omega\cdot(x_1-x_2)}\big(\omega\cdot(x_1-x_2)\big)^{M-j}(r^{-j}D^{2j}_\omega a(\omega))d\omega+O(\frac{\|D^{2(M+1)}_\omega a\|_{L^{\infty}}}{r^{M+1}}),
\end{multline*}
where $D^{2j}_{\omega}$ are differential operators of order $2j$ whose coefficients, as functions are polynomials in $x_1, x_2$ (i.e. a finite sum of $P_1(x_1)P_2(x_2)\psi(\omega)$).
\end{lem}
We shall see that the separation between $x_1, x_2$ in coefficients of $D^{2j}_{\omega}$ is crucial for the proof of Lemma \ref{lem-Mattila} below. 
\begin{proof}[Proof of Lemma \ref{lem-reverse-stationary-phase}]
 The classical stationary phase expansion (see, e.g.\cite{Ste93}\cite{Wol03}\cite{Sogge17}) states that, if the amplitude function $\psi\in C_0^\infty(S^{d-1})$ is supported on a local chart, then
 \begin{multline}\label{classical-stationary-phase}
 C_d\int e^{-2\pi ir\omega\cdot(x_1-x_2)}\psi(\omega)d\omega = (r|x_1-x_2|)^{-\frac{d-1}{2}}e^{-2\pi ir|x_1-x_2|}\\\sum_{j=0}^{M}(r|x_1-x_2|)^{-j} D^{2j}_{\omega}\psi(\omega_{x_1,x_2})+O(r^{-(M+1)}),
\end{multline}
where $D^{2j}_{\omega}$ are differential operators of order $2j$.

Now we consider a more complicated amplitude function $$\big(\omega\cdot(x_1-x_2)\big)^{M} a(\omega)\in C_0^\infty(S^{d-1}),$$ and move the term $j=0$ to the left hand side of the expansion. This gives
\begin{multline}\label{reverse-j=0}
e^{-2\pi i r|x_1-x_2|}|x_1-x_2|^{M-\frac{d-1}{2}}a(\omega_{x_1,x_2})\\=C_d\,r^{\frac{d-1}{2}}\int e^{-2\pi ir\omega\cdot(x_1-x_2)}\big(\omega\cdot(x_1-x_2)\big)^{M}a(\omega)d\omega-e^{-2\pi ir|x_1-x_2|}\\\sum_{j=1}^{M}r^{-j}|x_1-x_2|^{M-\frac{d-1}{2}-j} D^{2j}_{\omega}a(\omega_{x_1,x_2})+O(r^{-M-1}\|D^{2(M+1)}_\omega a\|_{L^{\infty}}),
\end{multline}
where coefficients of $D^{2j}_{\omega}$ are polynomials in $x_1, x_2$ because $\big(\omega\cdot(x_1-x_2)\big)^{M}$ is a polynomial in $x_1, x_2$. Here and in the rest of the proof, $D^{2j}_{\omega}$ may vary from line to line, but their coefficients are always polynomials in $x_1, x_2$.

By repeating the process \eqref{reverse-j=0} with 
$$\big(\omega\cdot(x_1-x_2)\big)^{M-1} D^2_\omega a(\omega)\in C_0^\infty(S^{d-1})$$
as the amplitude function, we obtain
\begin{multline}\label{reverse-j=1}
r^{-1}e^{-2\pi i r|x_1-x_2|}|x_1-x_2|^{M-\frac{d-1}{2}-1}D^2_\omega a(\omega_{x_1,x_2})\\=C_d\left(r^{\frac{d-1}{2}-1}\int e^{-2\pi ir\omega\cdot(x_1-x_2)}\big(\omega\cdot(x_1-x_2)\big)^{M-1}D^2_\omega a(\omega)d\omega\right)-e^{-2\pi ir|x_1-x_2|}\\\sum_{j=2}^{M}r^{-j}|x_1-x_2|^{M-\frac{d-1}{2}-j} D^{2j}_{\omega}a(\omega_{x_1,x_2})+O(r^{-M-1}\|D^{2(M+1)}_\omega a\|_{L^{\infty}}).
\end{multline}
By plugging \eqref{reverse-j=1} into \eqref{reverse-j=0}, we obtain
\begin{multline}\label{reverse-j=3}
e^{-2\pi i r|x_1-x_2|}|x_1-x_2|^{M-\frac{d-1}{2}}a(\omega_{x_1,x_2})\\=C_d\left(\sum_{j=0,1}r^{\frac{d-1}{2}-j}\int e^{-2\pi ir\omega\cdot(x_1-x_2)}\big(\omega\cdot(x_1-x_2)\big)^{M-j}D^{2j}_\omega a(\omega)d\omega\right)-e^{-2\pi ir|x_1-x_2|}\\\sum_{j=2}^{M}r^{-j}|x_1-x_2|^{M-\frac{d-1}{2}-j} D^{2j}_{\omega}a(\omega_{x_1,x_2})+O(r^{-M-1}\|D^{2(M+1)}_\omega a\|_{L^{\infty}}).
\end{multline}

Then we continue repeating this process with
$$\big(\omega\cdot(x_1-x_2)\big)^{M-j} D^{2j}_\omega a(\omega)\in C_0^\infty(S^{d-1}), j=1,\dots, M,$$
as amplitude functions, and eventually obtain
\begin{multline*}
e^{-2\pi i r|x_1-x_2|}|x_1-x_2|^{M-\frac{d-1}{2}}a(\omega_{x_1,x_2})\\=C_d\,\sum_{j=0}^M r^{\frac{d-1}{2}-j}\int e^{-2\pi ir\omega\cdot(x_1-x_2)}\big(\omega\cdot(x_1-x_2)\big)^{M-j}D^{2j}_\omega a(\omega)d\omega+O(r^{-M-1}\|D^{2(M+1)}_\omega a\|_{L^{\infty}})\\
= C_d\,\sum_{j=0}^M r^{\frac{d-1}{2}}\int e^{-2\pi ir\omega\cdot(x_1-x_2)}\big(\omega\cdot(x_1-x_2)\big)^{M-j}(r^{-j}D^{2j}_\omega a(\omega))d\omega+O(r^{-M-1}\|D^{2(M+1)}_\omega a\|_{L^{\infty}}),	
\end{multline*}
that completes the proof of Lemma \ref{lem-reverse-stationary-phase}. 
\end{proof}

\begin{proof}[Proof of Lemma \ref{lem-Mattila}]
By definition of $\nu_k$ in \eqref{def-nu-k},
\begin{equation}\label{hat-nu-k}\widehat{\nu_k}(r)=\iint e^{-2\pi i r|x_1-x_2|}|x_1-x_2|^{M-\frac{d-1}{2}}b_k(x_1,\omega_{x_1,x_2})b_k(x_2,\omega_{x_1,x_2})d\mu_1(x_1)d\mu_2(x_2).\end{equation}
By Lemma \ref{lem-reverse-stationary-phase}, the assumption $r\approx 2^{Tk}$ and \eqref{derivative-b}, there exists $M=M_{s,d}$ such that
$$e^{-2\pi i r|x_1-x_2|}|x_1-x_2|^{M-\frac{d-1}{2}}b_k(x_1, \omega_{x_1,x_2})b_k(x_2, \omega_{x_1,x_2})$$ 
$$=C_d\,\sum_{j=0}^M r^{\frac{d-1}{2}}\int e^{-2\pi ir\omega\cdot(x_1-x_2)}\big(\omega\cdot(x_1-x_2)\big)^{M-j}\Big(r^{-j}D^{2j}_\omega \big(b_k(x_1, \omega)b_k(x_2, \omega)\big)\Big)d\omega+O(r^{-d-1}).$$
Since coefficients of $D^{2j}_{\omega}$ are polynomials in $x_1, x_2$, this is a finite sum of
$$\begin{aligned}&r^{\frac{d-1}{2}}\int e^{-2\pi ir\omega\cdot(x_1-x_2)}P_1(x_1)P_2(x_2)\big(r^{-j_1}\partial_\omega^{2j_1} b_k(x_1, \omega)\big)\big(r^{-j_2}\partial_\omega^{2j_2} b_k(x_2, \omega)\big)\psi(\omega)d\omega\\=&r^{\frac{d-1}{2}}\int \Big(e^{-2\pi ir\omega\cdot x_1}P_1(x_1)\big(r^{-j_1}\partial_\omega^{2j_1} b_k(x_1, \omega)\big)\Big)\Big(e^{2\pi ir\omega\cdot x_2}P_2(x_2)\big(r^{-j_2}\partial_\omega^{2j_2} b_k(x_2, \omega)\big)\Big)\psi(\omega)d\omega,\end{aligned}$$
where $j=j_1+j_2$. Thanks to \eqref{derivative-b}, it suffices to consider the term $j=0$ when $r\approx 2^{Tk}$. Without loss of generality, we may also omit $P_1(x_1), P_2(x_2)$. Then, by Cauchy-Schwarz,
$$\begin{aligned}|\widehat{\nu_k}(r)|^2\lesssim & r^{d-1}\left|\int_{S^{d-1}} \Big(\int e^{-2\pi ir\omega\cdot x} b_k(x, \omega)d\mu_1(x)\Big)\Big(\int e^{2\pi ir\omega\cdot x} b_k(x, \omega)d\mu_2(x)\Big)d\omega\right|^2\\\leq &r^{d-1} \int \left|\int e^{-2\pi irx\cdot\omega}b_k(x,\omega)\,d\mu_1(x)\right|^2d\omega\cdot \int \left|\int e^{-2\pi irx\cdot\omega}b_k(x,\omega)\,d\mu_2(x)\right|^2d\omega\end{aligned}.$$
Hence, by one more Cauchy-Schwarz,
$$\int_{r\approx 2^{Tk}}|\widehat{\nu_k}(r)|^2dr\lesssim \prod_{j=1,2}\left(\int_{r\approx 2^{Tk}}\left(\int_{S^{d-1}} \left|\int e^{-2\pi irx\cdot\omega}b_k(x,\omega)\,d\mu_j(x)\right|^2d\omega\right)^2r^{d-1}dr\right)^{1/2},$$
that completes the proof of Lemma \ref{lem-Mattila}.

\end{proof}

\section{Cut the branching interval from the left}\label{sec-left-to-right}
As $\mu_1,\mu_2$ have been split in Lemma \ref{lem-Mattila}, from now we only consider $\mathcal{M}_k(\mu)$. More generally, for $0\leq l_1\leq h<l_2\leq k$, define
\begin{equation}\label{def-M-k-h}\mathcal{M}^{[l_1,l_2]}(\mu^{\D_h}):=\sup_{Q\in\D_h(\mu)}\int_{r\approx 2^{T(l_2-h)}}\left(\int \left|\int e^{-2\pi irx\cdot\omega}b^{[l_1,l_2]}_h(x,\omega)\,d\mu^Q(x)\right|^2d\omega\right)^2r^{d-1}dr.\end{equation}
As a remark, the branching interval for $\mathcal{M}^{[l_1,l_2]}(\mu^{\D_h})$ is $[h,l_2]$, but $l_1<h$ can still make a difference through the microlocal cutoff $b^{[l_1,l_2]}_h$.
\begin{lem}\label{lem-left-to-right} 
	\begin{equation}\label{estimate-initial}\M_k(\mu)\lessapprox 2^{-F_\mu(k_0)} \M^{[0,k]}(\mu^{\D_{k_0}}).\end{equation}	
Moreover, if there exists $[h_1,h_2]\subset[l_1,l_2]$ with $j_{[h_1,h_2]}\in(h, h+d_{[h_1,h_2]}]$ and $j_{[h_1,h_2]}+d^{[l_1,l_2]}_h\le l_2$. Then
	\begin{equation}\label{estimate-left-to-right}\mathcal{M}^{[l_1,l_2]}(\mu^{\D_h})\lessapprox 2^{2F_\mu(h)-F_\mu(j_{[h_1,h_2]})-\min_{[h,j_{[h_1,h_2]}]}F_\mu} \mathcal{M}^{[l_1,l_2]}(\mu^{\D_{j_{[h_1,h_2]}}}).\end{equation}	
\end{lem}
Notice that the branching interval is cut from $[0,k]$ to $[k_0, k]$ in \eqref{estimate-initial}, and cut from $[h,l_2]$ to $[j_{[h_1,h_2]},l_2]$ in \eqref{estimate-left-to-right}.

The role of \eqref{estimate-initial} is to provide an absolute power-saving as an initial step in the reduction in Section \ref{sec-proof-main}, due to $F_\mu(j)\geq (s-d/2)Tj$. The estimate \eqref{estimate-left-to-right} will be used in subsequent steps for either power-saving or power-loss.

\begin{proof}
    [Proof of Lemma \ref{lem-left-to-right}]
We only prove \eqref{estimate-initial} because the argument for \eqref{estimate-left-to-right} is identical.
  
Write out all squares in $\mathcal{M}_{k}(\mu)$ defined in \eqref{def-M} to obtain
$$\int\left(\iiint
	e^{-2\pi ir(\omega_1\cdot(x_1-x_1')+\omega_2\cdot(x_2-x_2'))}
	\prod_{j=1,2}b_k(x_j,\omega_j)b_k(x_j',\omega_j) d\omega_1 d\omega_2\,\phi(\frac{r}{2^{Tk}})r^{d-1}dr\right)(d\mu)^{4}.$$
Decompose
$$\mu=\sum_{Q\in\mathcal{D}_{k_0}(\mu)}\mu(Q)\mu_Q$$
and consider $(x_1, x_1',x_2, x_2')\in Q_1\times Q_1'\times Q_2\times Q_2'$. As $r\approx 2^{Tk}$ and $b_k$ satisfies \eqref{derivative-b} with $k_0\leq k/3$, integrations by parts shows that one can only consider quadruples $(Q_1, Q_1', Q_2, Q_2')$ whose centers satisfy
\begin{equation}\label{relation}\begin{aligned}|\pi_{\omega_1^\perp}(c_{Q_1}-c_{Q_1'})|,\,|\pi_{\omega_2^\perp}(c_{Q_2}-c_{Q_2'})|&\lesssim 2^{-Tk_0},\\|\omega_1\cdot(c_{Q_1}-c_{Q_1'})
       +\omega_2\cdot(c_{Q_2}-c_{Q_2'})|&\lesssim 2^{-Tk_0},\end{aligned}\end{equation}
       where $\pi_{\omega^\perp}$ refers to the orthogonal projection to $\omega^\perp$. This means, if we fix $\omega_1, \omega_2, r, Q_1, Q_2$, then $Q_1'$ must lie in a $2^{-Tk_0}\times\cdots\times 2^{-Tk_0}\times 1$ tube $T_{k_0}$ of direction $\omega_1$. Notice that, due to the cutoff $\Theta_{Q, 0}$ in \eqref{def-b-k-initial}, this tube $T_{k_0}$ must be good. If, in addition, $Q_1'$ is also fixed, then $Q_2'$ is fixed as well. Denote this relation by
$$(Q_1, Q_2)\sim (Q_1', Q_2').$$
Also recall that the regularization of $\mu$ guarantees that all cubes $Q\in\mathcal{D}_l(\mu)$ have approximately the same measure. Therefore, up to a negligible error, $\mathcal{M}_k(\mu)$ equals
\begin{multline*}
	\int_{r\approx 2^{Tk}}\iint\sum_{\substack{Q_1, Q_2, Q_1', Q_2'\in\mathcal{D}_{k_0}(\mu)\\(Q_1, Q_2)\sim (Q_1', Q_2')}}\mu(Q)^4 \\\left(\prod_{j=1,2}\int e^{-2\pi irx\cdot\omega_j}b_k(x,\omega_j)\,d\mu_{Q_j}(x)\int e^{2\pi irx\cdot\omega_j}b_k(x,\omega_j)\,d\mu_{Q_j'}(x)\right)d\omega_1 d\omega_2r^{d-1}dr.
\end{multline*}
By the inequality $2ab\leq a^2+b^2$, this integral is
\begin{multline*}
	\leq \int_{r\approx 2^{Tk}}\iint\sum_{\substack{Q_1, Q_2, Q_1', Q_2'\in\mathcal{D}_{k_0}(\mu)\\(Q_1, Q_2)\sim (Q_1', Q_2')}}\mu(Q)^4\Big(\left|\prod_{j=1,2}\int e^{-2\pi irx\cdot\omega_j}b_k(x,\omega_j)\,d\mu_{Q_j}(x)\right|^2\\+\left|\prod_{j=1,2}\int e^{-2\pi irx\cdot\omega_j}b_k(x,\omega_j)\,d\mu_{Q_j'}(x)\right|^2\Big) d\omega_1 d\omega_2\,r^{d-1}dr
\end{multline*}
\begin{multline*}
	\lesssim \int_{r\approx 2^{Tk}}\iint\sum_{Q_1, Q_2\in\mathcal{D}_{k_0}(\mu)}\mu(Q)^4\cdot \#\{(Q_1', Q_2'): (Q_1', Q_2')\sim (Q_1, Q_2)\}\\\left|\prod_{j=1,2}\int e^{-2\pi irx\cdot\omega_j}b_k(x,\omega_j)\,d\mu_{Q_j}(x)\right|^2d\omega_1 d\omega_2r^{d-1}dr,
\end{multline*}
which is, by our counting on $Q_1', Q_2'$ and Cauchy-Schwarz,
$$\lesssim\sup_{T_{k_0}\text{ good}} \#\{Q\subset T_{k_0}\} \int_{r\approx 2^{Tk}}\left(\sum_{Q\in\mathcal{D}_{k_0}(\mu)}\mu(Q)^2\int\left|\int e^{-2\pi irx\cdot\omega}b_k(x,\omega)\,d\mu_Q(x)\right|^2d\omega\right)^2 r^{d-1}dr$$
$$\lessapprox \sup_{T_{k_0}\text{ good}} \mu(T_{k_0})\cdot \sup_{Q\in \D_{k_0}(\mu)}\mu(Q)\int_{r\approx 2^{Tk}}\left(\int\left|\int e^{-2\pi irx\cdot\omega}b_k(x,\omega)\,d\mu_Q(x)\right|^2d\omega\right)^2 r^{d-1}dr.$$
By \eqref{good tube}, \eqref{measure-of-cube}, and rescaling, it is
$$\lessapprox2^{-F_\mu(k_0)-\min_{[0,k_0]}F_\mu} \sup_{Q\in \D_{k_0}(\mu)}\int_{r\approx 2^{T(k-k_0)}}\left(\int\left|\int e^{-2\pi irx\cdot\omega}b_k(2^{-Tk_0}x,\omega)\,d\mu^Q(x)\right|^2d\omega\right)^2 r^{d-1}dr.$$
It remains to replace $b_k$ by $b^{[0,k]}_{k_0}$. By the observation (ii) in Section \ref{sec-distance-measure}, the first factor in \eqref{def-b-k-initial} and all terms for $j\leq k_0$ in \eqref{def-b-k} are functions in $\omega$ independent in $x$. Therefore we can replace $b_k$ by $b^{[0,k]}_{k_0}$ as an upper bound to complete the proof of \eqref{estimate-initial}. We refer to \cite{Liu26+} for a similar idea of this last step on the pinned distance problem.

For \eqref{estimate-left-to-right}, the proof is identical. The condition $j_{[h_1,h_2]}+d^{[l_1,l_2]}_h\le l_2$ guarantees the integration by parts and the condition $j_{[h_1,h_2]}\in(h, h+d_{[h_1,h_2]}]$ ensures all $2^{-Tj_{[h_1,h_2]}}\times\cdots\times 2^{-Tj_{[h_1,h_2]}}\times 2^{-Th}$ tubes under consideration are good. We leave details to readers.
\end{proof}

\section{Cut the branching interval from the right}\label{sec-right-to-left}
Discussions in previous sections are on
$$\int_{r\approx 2^{Tk}}|\widehat{\nu_{k}}(r)|^2dr.$$
In this section, we consider a related quantity. Given $k'\leq k-l$ and $Q_1, Q_2\in\D_l(\mu)$ with $\dist(Q_1, Q_2)\approx 2^{-Tl}$, similar to \eqref{def-nu-k}, one can define a distance measure $\nu_{k'}^{Q_1, Q_2}$ by
$$\iint f(|x_1-x_2|)|x_1-x_2|^{M_{d,s}-\frac{d-1}{2}}b_l^{[l,l+k']}(x_1,\omega_{x_1,x_2})b_l^{[l,l+k']}(x_2,\omega_{x_1,x_2})d\mu^{Q_1}(x_1)d\mu^{Q_2}(x_2).$$
Let
\begin{equation}\label{def-E} \mathcal{E}_{k'}(\mu^{\D_l}):=\sup_{\substack{Q_1, Q_2\in\D_l(\mu)\\\dist(Q_1, Q_2)\approx 2^{-Tl}}}\int_{r\leq 2^{Tk'}}|\widehat{\nu_{k'}^{Q_1, Q_2}}(r)|^2dr.\end{equation}
An equivalent definition of $\mathcal{E}_{k'}(\mu^{\D_l})$ is the following. Up to a negligible error,
$$\mathcal{E}_{k'}(\mu^{\D_l})\approx \int|\widehat{\nu_{k'}^{Q_1,Q_2}}(r)|^2\psi(r/2^{Tk'})dr=2^{Tk'} \iint \hat{\psi}(2^{Tk'}(t_1-t_2))d\nu_{k'}^{Q_1,Q_2}(t_1)d\nu_{k'}^{Q_1,Q_2}(t_2),$$
where $\psi, \hat{\psi}$ are nonnegative smooth functions and $\hat{\psi}$ has bounded support. Therefore,
\begin{multline}
    \label{def-E-equi}
       \mathcal{E}_{k'}(\mu^{\D_l}) \approx \sup_{\substack{Q_1, Q_2\in\D_l(\mu)\\\dist(Q_1, Q_2)\approx 2^{-Tl}}}2^{Tk'}\iiiint\limits_{||x_1-x_2|-|x_1'-x_2'||\lesssim 2^{-Tk'}}b^{[l,l+k']}_l(x_1,\omega_{x_1,x_2})\cdots b^{[l,l+k']}_l(x_2',\omega_{x_1',x_2'})\\d\mu^{Q_1}(x_1)d\mu^{Q_2}(x_2)d\mu^{Q_1}(x_1')d\mu^{Q_2}(x_2').	
\end{multline}
\begin{lem}\label{lem-proj}
Suppose $\mu$ is $2^{-T(l+k')}$-regularized and $j:=j_{[l,l+k']}\geq l+k'/2$. Then
	\begin{equation}\mathcal{E}_{k'}(\mu^{\D_l})\lessapprox 2^{-\min_{[j,l+k']}(F_\mu-F_\mu(j))}\cdot \mathcal{E}_{j-l}(\mu^{\D_l}).\end{equation}	
\end{lem}
Notice that the branching interval is cut from $[l,l+k']$ to $[l,j_{[l,l+k']}]$.
\begin{proof}[Proof of Lemma \ref{lem-proj}] Consider the equivalent definition of $\mathcal{E}_{k'}(\mu^{\D_l})$ in \eqref{def-E-equi} with $Q_1,Q_2\in\D_l(\mu)$ as in \eqref{def-E}. Notice that, for fixed $(x_1', x_2')\in\supp\mu^{Q_1}\times\supp\mu^{Q_2}$ with $t:=|x_1'-x_2'|\approx 1$, one only needs to consider
$$(x_1, x_2)\in A^{x_1', x_2'}_{k'}:=\{(x_1, x_2): |x_1-x_2|\in (t-2^{-Tk'}, t+2^{-Tk'})\}.$$
Lemma \ref{lem-proj} basically states that one can thicken $A^{x_1', x_2'}_{k'}$ to $A^{x_1', x_2'}_{j-l}$ with reasonable cost.

Decompose the unit cubes into $2^{-T(j-l)}$-cubes. Since $j\geq l+k'/2$, the intersection between a $2^{-Tk'}$-annulus and a $2^{-T(j-l)}$-cube is approximately a $2^{-T(j-l)}\times\cdots\times 2^{-T(j-l)}\times 2^{-Tk'}$ slab. We refer to \cite{Orp17}\cite{KS18} for previous work on this geometric observation.

Therefore, for each $X_1\times X_2\in\D_{j-l}(\mu^{Q_1})\times\D_{j-l}(\mu^{Q_2})$ that intersects $A^{x_1', x_2'}_{k'}$, one can cover $(X_1\times X_2)\cap A^{x_1', x_2'}_{k'}$ by a union of products of slabs. More precisely,
\begin{equation}\label{sup-b}
\begin{aligned}
	&\iint_{A^{x_1', x_2'}_{k'}}b_l^{[l,l+k']}(x_1,\omega_{x_1,x_2})b_l^{[l,l+k']}(x_1,\omega_{x_1,x_2})d\mu^{Q_1}_{X_1}(x_1)d\mu^{Q_2}_{X_2}(x_2)\\\lesssim &\sup_{(x_1,x_2)\in X_1\times X_2}b_l^{[l,l+k']}(x_1,\omega_{x_1,x_2})\,b_l^{[l,l+k']}(x_2,\omega_{x_1,x_2})\cdot\sum_{S_\omega} \mu^{Q_1}_{X_1}(S_\omega)\cdot\mu^{Q_2}_{X_2}(S_\omega+t\omega),
\end{aligned}
\end{equation}
where $\omega:=\omega_{x_2,x_1}$ for some fixed $(x_1, x_2)\in (X_1\times X_2)\cap A^{x_1', x_2'}_{k'}$, and the sum is taken over almost disjoint $2^{-T(j-l)}\times\cdots\times 2^{-T(j-l)}\times 2^{-Tk'}$ slabs in $X_1$ perpendicular to $\omega$, which are good on average by our microlocal cutoff. Notice that the argument here does not depend on the choice of $(x_1, x_2)$ for $\omega$.

The second factor in \eqref{sup-b} is, by Cauchy-Schwarz and \eqref{good slab},
$$\leq \left(\sum_{S_\omega} \mu^{Q_1}_{X_1}(S_\omega)^2\right)^{1/2}\left(\sum_{S_\omega} \mu^{Q_2}_{X_2}(S_\omega+t\omega)^2\right)^{1/2}\approx \prod_{i=1,2}\left(\int \mu^{Q_i}_{X_i}(S_\omega(x))d\mu^{Q_i}_{X_i}(x)\right)^{1/2} $$
$$\lessapprox  2^{-Tk'}\cdot 2^{-T(l-j)-\min_{[j,l+k']} (F_\mu-F_\mu(j))}.$$
Since $\mathcal{E}_{k'}(\mu^{\D_l})$ is a sum of \eqref{sup-b} over $X_1\times X_2$, the problem has been reduced to estimating
\begin{equation}\label{sum-sup}
	\sum_{(X_1\times X_2)\cap A^{x_1', x_2'}_{k'}\neq\emptyset}\mu^{Q_1}(X_1)\mu^{Q_2}(X_2)\sup_{(x_1,x_2)\in X_1\times X_2}b_l^{[l,l+k']}(x_1,\omega_{x_1,x_2})\,b_l^{[l,l+k']}(x_2,\omega_{x_1,x_2}).
\end{equation}
By observations (i)(ii) from Section \ref{sec-distance-measure}, we have $b^{[l,l+k']}_l\leq b^{[l,j]}_l$ that is a function on $\omega$ independent in $x$ on each $2^{-T(j-l)}$-cube $X$. Also, when $X_1\times X_2$ is fixed, $\{\omega_{x_1, x_2}: (x_1, x_2)\in X_1\times X_2\}$ is contained in a $2^{-T(j-l)}$-cap. Therefore, \eqref{sum-sup} is
\begin{equation}\label{int-sup}
\begin{aligned}
	&\lesssim \sum_{(X_1\times X_2)\cap A^{x_1', x_2'}_{k'}\neq\emptyset}\int_{X_1}\int_{X_2} b_l^{[l,j]}(x_1,\omega_{x_1,x_2})\,b_l^{[l,j]}(x_2,\omega_{x_1,x_2})\,d\mu^{Q_1}(x_1)\mu^{Q_2}(x_2)\\&\leq\iint_{A^{x_1', x_2'}_{j-l}} b_l^{[l,j]}(x_1,\omega_{x_1,x_2})\,b_l^{[l,j]}(x_2,\omega_{x_1,x_2})\,d\mu^{Q_1}(x_1)\mu^{Q_2}(x_2).
\end{aligned}
\end{equation}
Here we abuse the notation a little bit, because to control the supremum, one needs to enlarge the support of $\Theta$ from its original $2\cdot 2^{-(1-\epsilon)T(j-l)}$ (recall \eqref{def-Theta}) to $4\cdot 2^{-(1-\epsilon)T(j-l)}$ to define the microlocal cutoff in \eqref{int-sup}. But it makes no difference as the reduction terminates within finitely many steps independent in $k$.

Overall, 
\begin{multline*}
\mathcal{E}_{k'}(\mu^{\D_l})\lessapprox 2^{-\min_{[j,l+k']} (F_\mu-F_\mu(j))}\cdot 2^{T(j-l)}\iiiint\limits_{||x_1-x_2|-|x_1'-x_2'||\lesssim 2^{-T(j-l)}}\\b^{[l,j]}_l(x_1,\omega_{x_1,x_2})
\cdots b^{[l,j]}_l(x_2',\omega_{x_1',x_2'})\,d\mu^{Q_1}(x_1)\cdots d\mu^{Q_2}(x_2'),
\end{multline*}
which completes the proof by the equivalent definition of $\mathcal{E}_{j-l}(\mu^{\D_l})$ from \eqref{def-E-equi}.
\end{proof}

\section{Transitions between $\M$ and $\E$}\label{sec-transition}
\subsection{From $\E$ to $\M$}
\begin{lem}\label{lem-high-low}Suppose $\mu$ is $2^{-T(k'+l)}$-regularized. Then, for all $\epsilon k\leq h\leq k'$,
\begin{equation} \mathcal{E}_{k'}(\mu^{\D_l})\lessapprox \mathcal{E}_{h}(\mu^{\D_l})+ \sup_{h<l'\leq k'}\M^{[l,l+l']}(\mu^{\D_l}).\end{equation}
\end{lem}
\begin{proof}[Proof of Lemma \ref{lem-high-low}]
By the definition of $\mathcal{E}_{k'}(\mu^{\D_l})$ in \eqref{def-E},
\begin{equation}\label{E-lemma-high-low}\begin{aligned}
	\mathcal{E}_{k'}(\mu^{\D_l})&=\int_{r\leq 2^{Tk'}}|\widehat{\nu_{k'}^{\D_l}}(r)|^2dr\\&=\int_{r\leq 2^{(1-\epsilon^2)Tk'}}|\widehat{\nu_{k'}^{\D_l}}(r)|^2dr+\int_{2^{(1-\epsilon^2)Tk'}\leq r\leq 2^{Tk'}}|\widehat{\nu_{k'}^{\D_l}}(r)|^2dr.
\end{aligned}\end{equation}
By the remark after Lemma \ref{lem-Mattila}, the high frequency term is $\lessapprox \M^{[l,l+k']}(\mu^{\D_l})$. For the low frequency term, similar to our derivation of the equivalent definition \eqref{def-E-equi} for $\mathcal{E}_{k'}(\mu^{\D_l})$,
\begin{multline*}
\int_{r\leq 2^{(1-\epsilon^2)Tk'}}|\widehat{\nu_{k'}^{\D_l}}(r)|^2dr\lesssim 2^{(1-\epsilon^2)Tk'}\iiiint\limits_{||x_1-x_2|-|x_1'-x_2'||\lesssim 2^{-(1-\epsilon^2)Tk'}}b^{[l,l+k']}_l(x_1,\omega_{x_1,x_2})\cdots \\\cdots b^{[l,l+k']}_l(x_2',\omega_{x_1',x_2'})d\mu^{Q_1}(x_1)\cdots d\mu^{Q_2}(x_2'),
\end{multline*}
where the implicit constant is independent in $k'$. Recall that $b^{[l,l+k']}_l\leq b^{[l,l+(1-\epsilon^2)k']}_l$ from Section \ref{sec-distance-measure}, so it becomes $\mathcal{E}_{(1-\epsilon^2)k'}(\mu^{\D_l})$. Repeating this process completes the proof.

\end{proof}

\subsection{From $\M$ to $\E$}

In this subsection we write $b_h$ for $b^{[l_1,l_2]}_h$ and $d_h$ for $d^{[l_1,l_2]}_h$ to simplify the notation, as the superscript will not change until the very end of the proof.
\begin{lem}\label{lem-from-M-to-E}
Suppose $\mu$ is $2^{-Tk}$-regularized, $0\leq l_1\leq h<l_2\leq k$, and $d_h\leq (l_2-h)/2$. Then
$$\mathcal{M}^{[l_1,l_2]}(\mu^{\D_h})\lessapprox 
2^{-2(F_\mu(l_2-2d_h)-F_\mu(h))}
\mathcal{M}^{[l_1,l_2]}(\mu^{\D_{l_2-2d_h}})+\max_{h<l\leq l_2-2d_h}
2^{-2(F_\mu(l)-F_\mu(h))}\mathcal{E}_{l_2-l}(\mu^{\D_l}).
$$
\end{lem}

\begin{proof}[Proof of Lemma \ref{lem-from-M-to-E}] 
By the Whitney decomposition, $\supp\mu^{\D_h}\times \supp\mu^{\D_h}$ equals
\begin{multline*}
	\{(X, X')\in \D_{l_2-h-2d_h}(\mu^{\D_h})\times \D_{l_2-h-2d_h}(\mu^{\D_h}): \dist(X,X')\lesssim 2^{-T(l_2-h-2d_h)}\}\cup \\\bigcup_{h\leq l< l_2-2d_h}\{(X, X')\in\D_{l-h}(\mu^{\D_h})\times \D_{l-h}(\mu^{\D_h}):\dist(X,X')\approx 2^{-T(l-h)}\}
\end{multline*}
$$:=\W_{l_2-h-2d_h}\cup \bigcup_{h\leq l< l_2-2d_h}\W_{l-h}.$$
Then $\mathcal{M}^{[l_1,l_2]}(\mu^{\D_h})$ defined in \eqref{def-M-k-h} can be written as
\begin{multline*}\int\Big|\sum_{h\leq l\leq l_2-2d_h}\sum_{(X, X')\in\W_{l-h}} \mu^{\D_h}(X)\mu^{\D_h}(X')\iiint e^{-2\pi ir\omega\cdot(x-x')}b_h(x,\omega)b_h(x',\omega) \\d\omega\,d\mu^{\D_h}_{X}(x)d\mu^{\D_h}_{X'}(x')\Big|^2 \psi(\frac{r}{2^{T(l_2-h)}})r^{d-1}dr,
	\end{multline*}
	which is, by Cauchy-Schwarz,
\begin{multline}\label{full-expansion-h1-h2}\lessapprox\sum_{h\leq l\leq l_2-2d_h}2^{-2(f_\mu(l)-f_\mu(h))}\max_{(X, X')\in\W_{l-h}}\int\Big| \iiint e^{-2\pi ir\omega\cdot(x-x')}b_h(x,\omega)b_h(x',\omega) \\d\omega\,d\mu^{\D_h}_{X}(x)d\mu^{\D_h}_{X'}(x')\Big|^2 \psi(\frac{r}{2^{T(l_2-h)}})r^{d-1}dr.
	\end{multline}
	
The diagonal term $l=l_2-2d_h$ in \eqref{full-expansion-h1-h2} is, after Cauchy-Schwarz and rescaling, 
\begin{multline*}
	 \lesssim 2^{-2(f_\mu(l)-f_\mu(h))}\cdot 2^{dT(l_2-h-2d_h)}\sup_{X\in\mathcal{D}_{l_2-h-2d_h}(\mu^{\D_h})}\int\\\left(\int \left|\int e^{-2\pi irx\cdot\omega}b_{l_2-2d_h}(x,\omega)\,d\mu^{2^{-Th}X}(x)\right|^2d\omega\right)^2\psi(\frac{r}{2^{2Td_h}})r^{d-1}dr
\end{multline*}
$$\lesssim 2^{-2 (F_\mu(l_2-2d_h)-F_\mu(h))}\cdot\M^{[l_1,l_2]}(\mu^{\D_{l_2-2d_h}}),$$
where $b_h$ can be replaced by $b_{l_2-2d_h}$ for the same reason as Section \ref{sec-left-to-right}.

For off-diagonal terms in \eqref{full-expansion-h1-h2}, thanks to \eqref{derivative-b}, we can apply the classical stationary phase \eqref{classical-stationary-phase} over $d\omega$, and only need to consider the first term of the expansion, that is
\begin{multline*}
	 2^{-2(f_\mu(l)-f_\mu(h))}\cdot 2^{(d-1)T(l-h)} \int\Big|\iint e^{\pm 2\pi ir|x-x'|}b_h(x,\omega_{x,x'})b_h(x',\omega_{x,x'})\\d\mu^{\D_h}_{X}(x)\,d\mu^{\D_h}_{X'}(x')\Big|^2\psi(\frac{r}{2^{T(l_2-h)}})dr
\end{multline*}
\begin{multline}\label{after-stationary-phase}\lesssim 2^{-2(F_\mu(l)-F_\mu(h))}\cdot 2^{T(l_2-l)}\cdot\iiiint|\hat{\psi}(2^{T(l_2-l)}
(|x_1-x_1'|\pm|x_2-x_2'|))|\\\prod_{j=1,2}b_h(2^{-T(l-h)}x_j,\omega_{x_j,x_j'})b_h(2^{-T(l-h)}x_j',\omega_{x_j,x_j'})\,d\mu^{2^{-Th}X_2'}(x_2')\cdots d\mu^{2^{-Th}X_1}(x_1),\end{multline}
where $(X_1, X_1'), (X_2, X_2')\in\W_{l-h}$.

We only	consider $|x_1-x_1'|-|x_2-x_2'|$ as $\hat{\psi}(2^{T(l_2-l)}(|x_1-x_1'|+|x_2-x_2'|))$ has fast decay. For the same reason at the end of Section \ref{sec-right-to-left}, the microlocal cutoff $b^{[l_1,l_2]}_h$ can be replaced by $b^{[l,l_2]}_l$. Therefore the equation \eqref{after-stationary-phase} is
$$\lesssim  2^{-2(F_\mu(l)-F_\mu(h))}\mathcal{E}_{l_2-l}(\mu^{\D_l})$$
by the equivalent definition of $\E_{l_2-l}$ in \eqref{def-E-equi}. This completes the proof of Lemma \ref{lem-from-M-to-E}.
\end{proof}

\section{Proof of the main theorem}\label{sec-proof-main}
Now we are ready to show
$$\int_{r\approx 2^{Tk}}|\widehat{\nu_k}(r)|^2dr\lesssim 2^{-\rho Tk}$$
to complete the proof of the main theorem. By Lemma \ref{lem-Mattila}, it suffices to show $\M_k(\mu)$ has power decay.

We first choose $k_0$. We claim that, there exists $\rho>0$ such that, for all large $k$, there exists $k_0\leq k/3$ such that
$$F_\mu(k_0)=\min_{[k_0,k]}F_\mu\geq \frac{(s-d/2)Tk}{4}-dT\geq \rho Tk.$$
To see this, as $F_\mu(j)\geq (s-d/2)Tj+error$, we have $F_\mu(j)> \frac{(s-d/2)Tk}{4}$ for all $j\geq \frac{k}{3}$. So
$$k_0:=\max \{j: F_\mu(j)\leq \frac{(s-d/2)Tk}{4}\}\leq k/3.$$
Since $F_\mu$ is $dT$-Lipschitz, it follows that
 $$F_\mu(k_0)\geq \frac{(s-d/2)Tk}{4}-dT\geq \rho Tk$$
 when $k$ is large enough in terms of $s,d$. Therefore, by Lemma \ref{lem-left-to-right},
\begin{equation}\label{eq-power-saving}\M_{k}(\mu)\lessapprox 2^{-F_\mu(k_0)}\M^{[0,k]}(\mu^{\D_{k_0}})\leq 2^{-\rho Tk}\M^{[0,k]}(\mu^{\D_{k_0}})\end{equation}
for an absolute $\rho=\rho_{s,d}>0$ independent in $k$. Since $F_\mu(k_0)=\min_{[k_0,k]}F_\mu$, it remains to show the following lemma, then choose $0<\delta\ll 1$, and finally take $0<\epsilon\ll\delta$.
\begin{lem}\label{lem-final}
Suppose $\mu$ is a $2^{-Tk}$-regularized measure. Then for $0\leq l_1\leq h<l_2\leq k$,
$$\mathcal M^{[l_1,l_2]}(\mu^{\mathcal D_h}), \mathcal E_{l_2-h}(\mu^{\mathcal D_h})\lesssim_{\epsilon, \delta} 2^{O(\epsilon/\delta)Tk+O(\delta\log(1/\delta))Tk}\cdot 2^{2F_\mu(h)-\min_{[h+\delta k,l_2]}F_\mu-\min_{[h,l_2-\delta k]}F_\mu}.$$

\end{lem}
\begin{proof}
	[Proof of Lemma \ref{lem-final}]
Notice that the branching intervals for $\mathcal M^{[l_1,l_2]}(\mu^{\mathcal D_h})$, $\mathcal E_{l_2-h}(\mu^{\mathcal D_h})$ are both $[h,l_2]$. The proof goes by reduction on its length. We shall focus on the main term first, leaving the discussion on accumulated error to the end. If $l_2-h\leq 10\delta k$, the reduction terminates by trivial estimates.
	
Recall $j_{[h_1, h_2]}, d_{[h_1, h_2]}, d^{[l_1,l_2]}_h, j^{[l_1,l_2]}_h$ from \eqref{def-j-h-1-h-2}-\eqref{def-j-max}.
	
{\bf Reduction on $\E$.} The discussion is on the location of $j:=j_{[h,l_2]}\in[h+\delta k,l_2-\delta k]$.

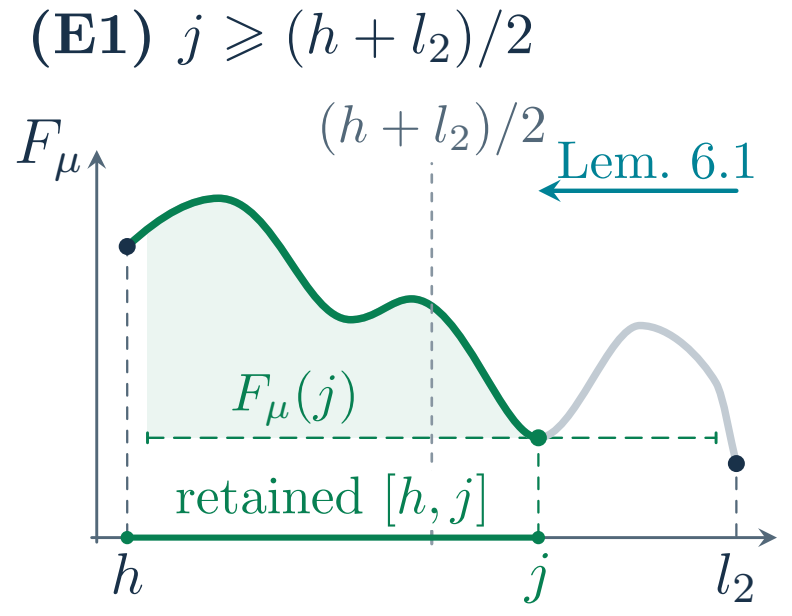
\begin{figure}[htbp]
\centering
\captionsetup{font=small,justification=raggedright,singlelinecheck=false,margin=0pt,skip=3pt}
\tikzset{energy panel/.style={x=.00650\linewidth,y=.55cm,>=stealth,line cap=round,line join=round,
  axis/.style={->,axiscolor,line width=.55pt},
  profile/.style={line width=1.4pt},
  projection/.style={dash pattern=on 2.2pt off 2.2pt,line width=.43pt},
  level/.style={dash pattern=on 3pt off 2pt,line width=.45pt},
  tick/.style={anchor=north,inner sep=1pt,text=ink,font=\small},
  annotation/.style={anchor=south,inner sep=1pt,font=\footnotesize},
  interval label/.style={anchor=south,fill=white,inner xsep=1pt,inner ysep=.5pt,font=\footnotesize},
  move arrow/.style={->,line width=.83pt}}}
  \definecolor{ink}{HTML}{193149}
  \definecolor{axiscolor}{HTML}{536879}
  \definecolor{contextcurve}{HTML}{C2CBD3}
  \definecolor{retain}{HTML}{078052}
  \definecolor{low}{HTML}{BE6909}
  \definecolor{move}{HTML}{008296}

\begin{minipage}[t]{.33\linewidth}
  \centering
  \begin{tikzpicture}[energy panel]
    \path[use as bounding box] (-4,-.95) rectangle (149,6.35);
    \node[anchor=west,inner sep=0pt,text=ink] at (0,6)
      {\small\bfseries (E1) \(j\geq (h+l_2)/2\)};
    \def\profilepath{(20,3.50)
      .. controls (27,3.90) and (33,4.08) .. (38,4.08)
      .. controls (49,4.08) and (55,2.62) .. (64,2.62)
      .. controls (69,2.62) and (72,2.87) .. (76,2.87)
      .. controls (86,2.87) and (94,1.20) .. (101,1.20)
      .. controls (108,1.20) and (115,2.55) .. (121,2.55)
      .. controls (127,2.55) and (133,2.17) .. (136,1.86)
      .. controls (137.5,1.70) and (138.8,1.18) .. (140,.89)}
    \begin{scope}
      \clip (24,1.20) rectangle (101,5.50);
      \fill[retain!8] \profilepath -- (140,1.20) -- (20,1.20) -- cycle;
    \end{scope}
    \draw[axis] (13,0) -- (148,0);
    \draw[axis] (14,-.055) -- (14,4.66);
    \node[anchor=east,inner sep=2pt,text=ink] at (13.8,4.66) {$F_\mu$};
    \foreach \x/\y/\col in
      {20/3.50/axiscolor,101/1.20/retain,140/.89/axiscolor}
      {\draw[projection,\col] (\x,.02) -- (\x,\y);}
    \draw[axiscolor!65,line width=.6pt] (80,-.075) -- (80,.075);
    \draw[retain,level] (24,1.20) -- (136,1.20);
    \foreach \x in {24,136}{\draw[retain,line width=.55pt] (\x,1.14) -- (\x,1.26);}
    \draw[profile,contextcurve] \profilepath;
    \begin{scope}
      \clip (20,-.02) rectangle (101,5.50);
      \draw[profile,retain] \profilepath;
    \end{scope}
    \foreach \x/\y/\col in
      {20/3.50/ink,101/1.20/retain,140/.89/ink}
      {\fill[\col] (\x,\y) circle[radius=1.6pt];}
    \draw[axiscolor!70,dash pattern=on 1.5pt off 2pt,line width=.5pt]
      (80,.02) -- (80,4.50);
    \node[font=\footnotesize,text=axiscolor,anchor=south,inner sep=1pt]
      at (80,4.52) {$(h+l_2)/2$};
    \node[tick] at (20,-.11) {$h$};
    \node[tick,text=retain] at (101,-.11) {$j$};
    \node[tick] at (140,-.11) {$l_2$};
    \node[annotation,text=retain] at (53,1.27) {$F_\mu(j)$};
    \draw[move arrow,move] (140,4.17) --
      node[pos=.40,annotation,above=1pt,text=move] {Lem.~\ref{lem-proj}} (101,4.17);
    \draw[retain,line width=1.15pt] (20,0) --
      node[midway,interval label,above=1pt,text=retain] {retained \([h,j]\)} (101,0);
    \foreach \x in {20,101} {\fill[retain] (\x,0) circle[radius=1.25pt];}
  \end{tikzpicture}
  \caption*{(E1): right cut by Lem. \ref{lem-proj}.}
\end{minipage}\hfill
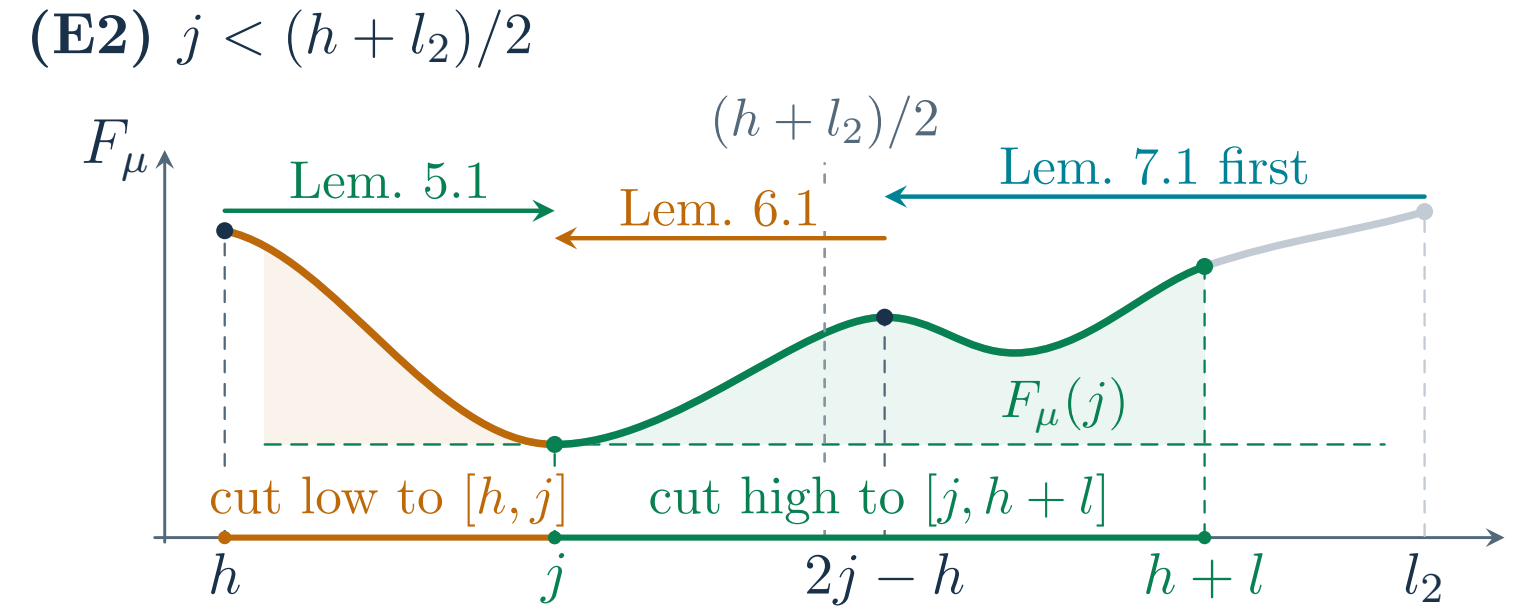
\begin{minipage}[t]{.65\linewidth}
  \centering
  \begin{tikzpicture}[energy panel]
    \path[use as bounding box] (-4,-.95) rectangle (149,6.35);
    \node[anchor=west,inner sep=0pt,text=ink] at (0,6)
      {\small\bfseries (E2)\ \(j<(h+l_2)/2\)};
    \def\profilepath{(20,3.69)
      .. controls (31,3.41) and (41,1.12) .. (53,1.12)
      .. controls (65,1.12) and (78,2.65) .. (86,2.65)
      .. controls (91,2.65) and (94,2.22) .. (99,2.22)
      .. controls (106,2.22) and (112,3.01) .. (118,3.26)
      .. controls (126,3.59) and (133,3.67) .. (140,3.92)}
    \begin{scope}
      \clip (24,1.12) rectangle (53,5.50);
      \fill[low!8] \profilepath -- (140,1.12) -- (20,1.12) -- cycle;
    \end{scope}
    \begin{scope}
      \clip (53,1.12) rectangle (118,5.50);
      \fill[retain!8] \profilepath -- (140,1.12) -- (20,1.12) -- cycle;
    \end{scope}
    \draw[axis] (13,0) -- (148,0);
    \draw[axis] (14,-.055) -- (14,4.66);
    \node[anchor=east,inner sep=2pt,text=ink] at (13.8,4.66) {$F_\mu$};
    \foreach \x/\y/\col in
      {20/3.69/axiscolor,53/1.12/retain,86/2.65/axiscolor,118/3.26/retain,140/3.92/contextcurve}
      {\draw[projection,\col] (\x,.02) -- (\x,\y);}
    \draw[retain,level] (24,1.12) -- (136,1.12);
    \draw[profile,contextcurve] \profilepath;
    \begin{scope}
      \clip (20,-.02) rectangle (53,5.50);
      \draw[profile,low] \profilepath;
    \end{scope}
    \begin{scope}
      \clip (53,-.02) rectangle (118,5.50);
      \draw[profile,retain] \profilepath;
    \end{scope}
    \foreach \x/\y/\col in
      {20/3.69/ink,53/1.12/retain,86/2.65/ink,118/3.26/retain,140/3.92/contextcurve}
      {\fill[\col] (\x,\y) circle[radius=1.6pt];}
    \draw[axiscolor!70,dash pattern=on 1.5pt off 2pt,line width=.5pt]
      (80,.02) -- (80,4.50);
    \node[font=\footnotesize,text=axiscolor,anchor=south,inner sep=1pt]
      at (80,4.6) {$(h+l_2)/2$};
    \node[tick] at (20,-.11) {$h$};
    \node[tick,text=retain] at (53,-.11) {$j$};
    \node[tick] at (86,-.11) {$2j-h$};
    \node[tick,text=retain] at (118,-.11) {$h+l$};
    \node[tick] at (140,-.11) {$l_2$};
    \node[annotation,text=retain] at (104,1.19) {$F_\mu(j)$};
    \begin{scope}[annotation/.append style={fill=white,inner xsep=1.5pt}]
      \draw[move arrow,move] (140,4.1) --
        node[midway,annotation,above=1pt,text=move] {Lem.~\ref{lem-high-low} first} (86,4.1);
      \draw[move arrow,low] (86,3.60) --
        node[midway,annotation,above=1pt,text=low] {Lem.~\ref{lem-proj}} (53,3.60);
      \draw[move arrow,retain] (20,3.93) --
        node[midway,annotation,above=1pt,text=retain] {Lem.~\ref{lem-left-to-right}} (53,3.93);
    \end{scope}
    \draw[low,line width=1.15pt] (20,0) --
      node[midway,interval label,above=1pt,text=low] {cut low to \([h,j]\)} (53,0);
    \foreach \x in {20,53} {\fill[low] (\x,0) circle[radius=1.25pt];}
    \draw[retain,line width=1.15pt] (53,0) --
      node[midway,interval label,above=1pt,text=retain] {cut high to \([j,h+l]\)} (118,0);
    \foreach \x in {53,118} {\fill[retain] (\x,0) circle[radius=1.25pt];}
  \end{tikzpicture}
  \caption*{(E2): split at $2j-h$ by Lem. \ref{lem-high-low} first, then cut low from right by Lem. \ref{lem-proj} and cut high from left by Lem. \ref{lem-left-to-right}.}
\end{minipage}
\end{figure}

{\bf Case (E1).} $j\geq(h+l_2)/2$.

In this case, we apply Lemma \ref{lem-proj} for
$$\E_{l_2-h}(\mu^{\D_{h}})\lessapprox 2^{-(\min_{[j,l_2]}F_\mu-F_\mu(j))}\E_{j-h}(\mu^{\D_{h}}).$$
By the definition of $j$ in \eqref{def-j-h-1-h-2}, the branching interval is shortened for  at least $\delta k$. So, if the main term in Lemma \ref{lem-final} has been verified on shorter branching intervals, then
$$\E_{l_2-h}(\mu^{\D_{h}})\lessapprox2^{-(\min_{[j,l_2]}F_\mu-F_\mu(j))}\cdot2^{2F_\mu(h)
-\min_{[h+\delta k,j]}F_\mu-\min_{[h,j-\delta k]}F_\mu}.$$
By the definition of $j=j_{[h,l_2]}$ in \eqref{def-j-h-1-h-2}, we have
$$F_\mu(j_{[h,l_2]})=\min_{[h+\delta k,l_2-\delta k]}F_\mu=\min_{[h+\delta k,j]}F_\mu,$$
so the exponent becomes
$$2F_\mu(h)-\min_{[j,l_2]}F_\mu-\min_{[h,j-\delta k]}F_\mu\leq 2F_\mu(h)-\min_{[h+\delta k,l_2]}F_\mu-\min_{[h,l_2-\delta k]}F_\mu,$$
as desired.

{\bf Case (E2).} $j<(h+l_2)/2$.

In this case, we apply Lemma \ref{lem-high-low} to obtain
\begin{equation}\label{combined-high-low}
\E_{l_2-h}(\mu^{\D_{h}})\lessapprox\E_{2(j-h)}(\mu^{\D_{h}})+\sup_{2(j-h)<l\leq l_2-h}\M^{[h,h+l]}(\mu^{\D_h}).
\end{equation}

For $\E_{2(j-h)}(\mu^{\D_{h}})$, if $2(j-h)<10\delta k$, nothing to prove; otherwise apply Lemma \ref{lem-proj} for
$$\E_{2(j-h)}(\mu^{\D_{h}})\lessapprox
2^{\,F_\mu(j)-\min_{[j,\,2j-h]}F_\mu}\E_{j-h}(\mu^{\D_{h}})$$
and the reduction follows in the same way as Case (E1).

For $\M^{[h,h+l]}(\mu^{\D_h})$, we may assume $l>10\delta k$. By the assumption $j<(h+l_2)/2$ and the trivial estimate $d^{[h,h+l]}_h\leq l/2$, we can apply Lemma \ref{lem-left-to-right} to obtain
$$\M^{[h,h+l]}(\mu^{\D_h})\lessapprox 2^{\,2F_\mu(h)-F_\mu(j)-\min_{[h,j]}F_\mu}
\M^{[h,h+l]}(\mu^{\D_{j}}).$$
The new branching interval is $[j,h+l]$, cut from $[h,h+l]$ by $j_{[h,l_2]}-h\geq \delta k$. So, if the main term in Lemma \ref{lem-final} has been verified on shorter branching intervals, then
$$\M^{[h,h+l]}(\mu^{\D_h})\lessapprox 2^{\,2F_\mu(h)-F_\mu(j)-\min_{[h,j]}F_\mu}\cdot2^{2F_\mu(j)-\min_{[j+\delta k,h+l]}F_\mu-\min_{[j,h+l-\delta k]}F_\mu}.$$
Since $h+\delta k\leq j_{[h,l_2]}\leq h+l-\delta k\leq l_2-\delta k$, we have
$$F_\mu(j_{[h,l_2]})=\min_{[h+\delta k,l_2-\delta k]}F_\mu=\min_{[j,h+l-\delta k]}F_\mu,$$
so the exponent becomes
$$2F_\mu(h)-\min_{[j+\delta k,h+l]}F_\mu-\min_{[h,j]}F_\mu\leq 2F_\mu(h)-\min_{[h+\delta k,l_2]}F_\mu-\min_{[h,l_2-\delta k]}F_\mu,$$
as desired. 

{\bf Reduction on $\M$.} The discussion is on $d_h:=d^{[l_1,l_2]}_h$ and $j_h:=j^{[l_1,l_2]}_h$. 
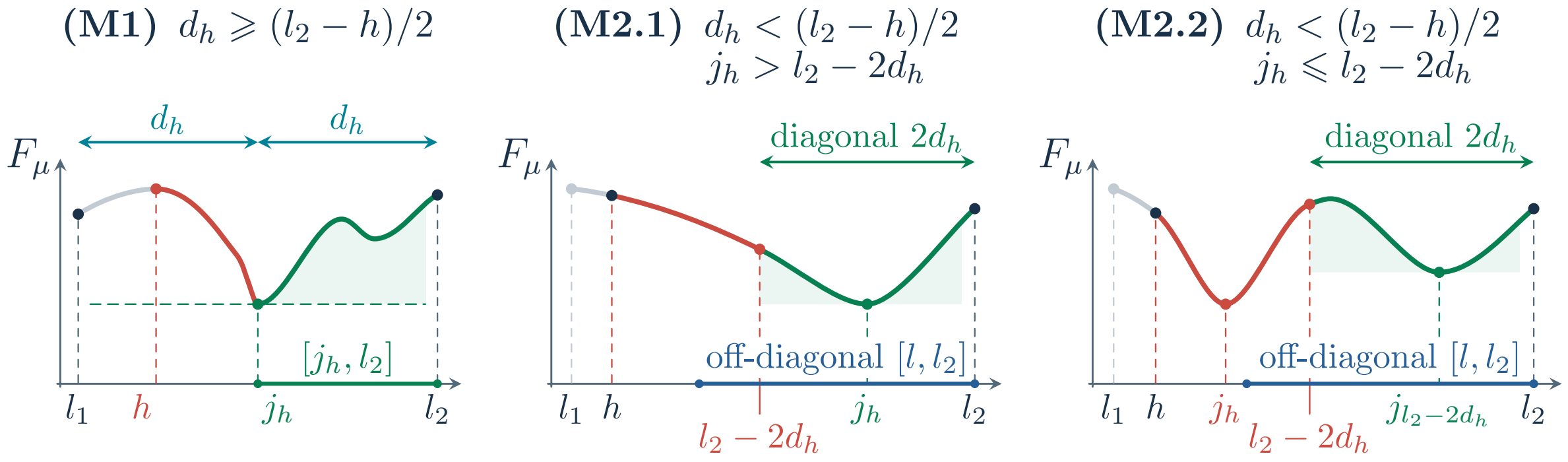
\begin{figure}[htbp]
\centering
\captionsetup{font=small,justification=raggedright,singlelinecheck=false,margin=0pt,skip=3pt}
\def\panelAxisTop{4.10}
\def\panelDimensionY{4.40}
\def\panelDiagonalY{4.04}
\tikzset{mattila panel/.style={x=.00650\linewidth,y=.55cm,baseline=0pt,>=stealth,line cap=round,line join=round,
  axis/.style={->,axiscolor,line width=.55pt},
  profile/.style={line width=1.4pt},
  projection/.style={dash pattern=on 2.2pt off 2.2pt,line width=.43pt},
  level/.style={dash pattern=on 3pt off 2pt,line width=.45pt},
  tick/.style={anchor=north,inner sep=1pt,text=ink,font=\footnotesize},
  annotation/.style={anchor=south,inner sep=1pt,font=\footnotesize},
  interval label/.style={anchor=south,fill=white,inner xsep=1pt,inner ysep=.5pt,font=\footnotesize},
  dimension/.style={<->,line width=.71pt},
  panel heading/.style={anchor=west,inner sep=0pt,text=ink,font=\small}}}
  \definecolor{ink}{HTML}{193149}
  \definecolor{axiscolor}{HTML}{536879}
  \definecolor{contextcurve}{HTML}{C2CBD3}
  \definecolor{retain}{HTML}{078052}
  \definecolor{direct}{HTML}{CB4B40}
  \definecolor{move}{HTML}{008296}
  \definecolor{offdiag}{HTML}{265F97}

\begin{minipage}[t]{.295\linewidth}
  \centering
  \begin{tikzpicture}[mattila panel]
    \path[use as bounding box] (-5,-1.40) rectangle (148,5.95);
    \node[panel heading] at (14,6.5) {\textbf{(M1)}\enspace \(d_h\geq(l_2-h)/2\)};
    \def\profilepath{(20,3.09)
      .. controls (30,3.43) and (40,3.55) .. (46,3.55)
      .. controls (58,3.55) and (67,2.75) .. (73,2.34)
      .. controls (75,2.20) and (76,1.96) .. (77,1.82)
      .. controls (78.10,1.65) and (79,1.45) .. (80,1.45)
      .. controls (89,1.45) and (99,3.00) .. (108,3.00)
      .. controls (112,3.00) and (115,2.64) .. (119,2.64)
      .. controls (127,2.64) and (134,3.28) .. (140,3.44)}
    \begin{scope}
      \clip (80,1.45) rectangle (136,5.50);
      \fill[retain!8] \profilepath -- (140,1.45) -- (20,1.45) -- cycle;
    \end{scope}
    \draw[axis] (13,0) -- (148,0);
    \draw[axis] (14,-.055) -- (14,\panelAxisTop);
    \node[anchor=east,inner sep=2pt,text=ink] at (13.8,\panelAxisTop) {$F_\mu$};
    \foreach \x/\y/\col in
      {20/3.09/axiscolor,46/3.55/direct,80/1.45/retain,140/3.44/axiscolor}
      {\draw[projection,\col] (\x,.02) -- (\x,\y);}
    \draw[retain,level] (24,1.45) -- (136,1.45);
    \draw[profile,contextcurve] \profilepath;
    \begin{scope}
      \clip (46,-.02) rectangle (80,5.50);
      \draw[profile,direct] \profilepath;
    \end{scope}
    \begin{scope}
      \clip (80,-.02) rectangle (140,5.50);
      \draw[profile,retain] \profilepath;
    \end{scope}
    \foreach \x/\y/\col in
      {20/3.09/ink,46/3.55/direct,80/1.45/retain,140/3.44/ink}
      {\fill[\col] (\x,\y) circle[radius=1.6pt];}
    \node[tick] at (20,-.11) {$l_1$};
    \node[tick,anchor=north east,text=direct] at (46,-.11) {$h$};
    \node[tick,anchor=north west,text=retain] at (80,-.11) {$j_h$};
    \node[tick] at (140,-.11) {$l_2$};
    \draw[dimension,move] (20,\panelDimensionY) --
      node[midway,annotation,above=1pt,text=move] {\(d_h\)} (80,\panelDimensionY);
    \draw[dimension,move] (80,\panelDimensionY) --
      node[midway,annotation,above=1pt,text=move] {\(d_h\)} (140,\panelDimensionY);
    \draw[retain,line width=1.15pt] (80,0) --
      node[midway,interval label,above=1pt,text=retain] {\([j_h,l_2]\)} (140,0);
    \foreach \x in {80,140} {\fill[retain] (\x,0) circle[radius=1.25pt];}
  \end{tikzpicture}
  \caption*{(M1): left cut by Lem.~\ref{lem-left-to-right}}
\end{minipage}\hfill
\begin{minipage}[t]{.69\linewidth}
\centering
\begin{minipage}[t]{.48\linewidth}
  \centering
  \begin{tikzpicture}[mattila panel]
    \path[use as bounding box] (-5,-1.40) rectangle (148,5.95);
    \node[panel heading] at (14,6.5) {\textbf{(M2.1)}\enspace \(d_h<(l_2-h)/2\)};
    \node[panel heading] at (60,5.76) {\(j_h>l_2-2d_h\)};
    \def\profilepath{(20,3.55)
      .. controls (24,3.53) and (29,3.47) .. (32,3.43)
      .. controls (48,3.21) and (61,2.90) .. (76,2.45)
      .. controls (87,2.12) and (99,1.45) .. (108,1.45)
      .. controls (117,1.45) and (132,2.76) .. (140,3.19)}
    \begin{scope}
      \clip (76,1.45) rectangle (136,5.50);
      \fill[retain!8] \profilepath -- (140,1.45) -- (20,1.45) -- cycle;
    \end{scope}
    \draw[axis] (13,0) -- (148,0);
    \draw[axis] (14,-.055) -- (14,\panelAxisTop);
    \node[anchor=east,inner sep=2pt,text=ink] at (13.8,\panelAxisTop) {$F_\mu$};
    \foreach \x/\y/\col in
      {20/3.55/contextcurve,32/3.43/direct,76/2.45/direct,108/1.45/retain,140/3.19/axiscolor}
      {\draw[projection,\col] (\x,.02) -- (\x,\y);}
    \draw[profile,contextcurve] \profilepath;
    \begin{scope}
      \clip (32,-.02) rectangle (76,5.50);
      \draw[profile,direct] \profilepath;
    \end{scope}
    \begin{scope}
      \clip (76,-.02) rectangle (140,5.50);
      \draw[profile,retain] \profilepath;
    \end{scope}
    \foreach \x/\y/\col in
      {20/3.55/contextcurve,32/3.43/ink,76/2.45/direct,108/1.45/retain,140/3.19/ink}
      {\fill[\col] (\x,\y) circle[radius=1.6pt];}
    \node[tick] at (20,-.11) {$l_1$};
    \node[tick] at (32,-.11) {$h$};
    \draw[direct,line width=.45pt] (76,0) -- (76,-.55);
    \node[tick,text=direct] at (76,-.64) {$l_2-2d_h$};
    \node[tick,text=retain] at (108,-.11) {$j_h$};
    \node[tick] at (140,-.11) {$l_2$};
    \draw[dimension,retain] (76,\panelDiagonalY) --
      node[midway,annotation,above=1pt,text=retain] {diagonal \(2d_h\)} (140,\panelDiagonalY);
    \draw[offdiag,line width=1.15pt] (58,0) --
      node[midway,interval label,above=1pt,text=offdiag] {off-diagonal \([l,l_2]\)} (140,0);
    \foreach \x in {58,140} {\fill[offdiag] (\x,0) circle[radius=1.25pt];}
  \end{tikzpicture}
\end{minipage}\hfill
\begin{minipage}[t]{.50\linewidth}
  \centering
  \begin{tikzpicture}[mattila panel]
    \path[use as bounding box] (-5,-1.40) rectangle (148,5.95);
    \node[panel heading] at (14,6.5) {\textbf{(M2.2)}\enspace \(d_h<(l_2-h)/2\)};
    \node[panel heading] at (60,5.76) {\(j_h\leq l_2-2d_h\)};
    \def\profilepath{(20,3.55)
      .. controls (24,3.46) and (28,3.31) .. (32,3.11)
      .. controls (39,2.76) and (46,1.45) .. (52,1.45)
      .. controls (59,1.45) and (68,3.09) .. (76,3.27)
      .. controls (78,3.32) and (80,3.37) .. (82,3.37)
      .. controls (92,3.37) and (104,2.03) .. (113,2.03)
      .. controls (123,2.03) and (133,2.82) .. (140,3.19)}
    \begin{scope}
      \clip (76,2.03) rectangle (136,5.50);
      \fill[retain!8] \profilepath -- (140,2.03) -- (20,2.03) -- cycle;
    \end{scope}
    \draw[axis] (13,0) -- (148,0);
    \draw[axis] (14,-.055) -- (14,\panelAxisTop);
    \node[anchor=east,inner sep=2pt,text=ink] at (13.8,\panelAxisTop) {$F_\mu$};
    \foreach \x/\y/\col in
      {20/3.55/contextcurve,32/3.11/direct,52/1.45/direct,76/3.27/direct,113/2.03/retain,140/3.19/axiscolor}
      {\draw[projection,\col] (\x,.02) -- (\x,\y);}
    \draw[profile,contextcurve] \profilepath;
    \begin{scope}
      \clip (32,-.02) rectangle (76,5.50);
      \draw[profile,direct] \profilepath;
    \end{scope}
    \begin{scope}
      \clip (76,-.02) rectangle (140,5.50);
      \draw[profile,retain] \profilepath;
    \end{scope}
    \foreach \x/\y/\col in
      {20/3.55/contextcurve,32/3.11/ink,52/1.45/direct,76/3.27/direct,113/2.03/retain,140/3.19/ink}
      {\fill[\col] (\x,\y) circle[radius=1.6pt];}
    \node[tick] at (20,-.11) {$l_1$};
    \node[tick] at (32,-.11) {$h$};
    \node[tick,text=direct] at (52,-.11) {$j_h$};
    \draw[direct,line width=.45pt] (76,0) -- (76,-.55);
    \node[tick,text=direct] at (76,-.64) {$l_2-2d_h$};
    \node[tick,text=retain] at (113,-.11) {$j_{l_2-2d_h}$};
    \node[tick] at (140,-.11) {$l_2$};
    \draw[dimension,retain] (76,\panelDiagonalY) --
      node[midway,annotation,above=1pt,text=retain] {diagonal \(2d_h\)} (140,\panelDiagonalY);
    \draw[offdiag,line width=1.15pt] (58,0) --
      node[midway,interval label,above=1pt,text=offdiag] {off-diagonal \([l,l_2]\)} (140,0);
    \foreach \x in {58,140} {\fill[offdiag] (\x,0) circle[radius=1.25pt];}
  \end{tikzpicture}
\end{minipage}
\par\caption*{(M2): split at $l_2-2d_h$ by Lem.~\ref{lem-from-M-to-E}, the off-diagonals then enter the $\E$ reduction; the diagonal is either cut from left by Lem.~\ref{lem-left-to-right} or spilt by Lem.~\ref{lem-from-M-to-E} again.}
\end{minipage}
\end{figure}

We shall need the following observation: if $j_{[h_1,h_2]}<j_{[h_1',h_2']}$ and $d_{[h_1,h_2]}\geq d_{[h_1',h_2']}$, then, by the minimality in \eqref{def-j-h-1-h-2}, there is no way for $j_{[h_1,h_2]}\in (j_{[h_1',h_2']}-d_{[h_1',h_2']}+\delta k, j_{[h_1',h_2']})$. In other wodrs, we must have
\begin{equation}\label{jh-jh'}j_{[h_1',h_2']}-j_{[h_1,h_2]}>d_{[h_1',h_2']}-\delta k.\end{equation}
In particular, if $j_{h'}^{[l_1', l_2']}>j_h^{[l_1, l_2]}$ with $[l_1', l_2']\subset [l_1, l_2]$ and $l_2'-h'>10\delta k$, then the maximality of $d_h^{[l_1, l_2]}$ in \eqref{def-d-max} guarantees that $d_h^{[l_1, l_2]}>d_{h'}^{[l_1', l_2']}$. So \eqref{jh-jh'} implies that
\begin{equation}\label{jh-jjh}j_{h'}^{[l_1', l_2']}-j_h^{[l_1, l_2]}>d_{h'}^{[l_1', l_2']}-\delta k.\end{equation}

{\bf Case (M1).} $d_h\geq (l_2-h)/2$.  

In this case, we have $h< j_h\leq l_2-d_h\leq h+d_h$, so Lemma \ref{lem-left-to-right} applies, giving
$$\M^{[l_1,l_2]}(\mu^{\D_h})\lessapprox 2^{\,2F_\mu(h)-F_\mu(j_h)-\min_{[h,j_h]}F_\mu}\M^{[l_1,l_2]}(\mu^{\D_{j_h}}).$$
The branching interval is cut from $[h, l_2]$ to $[j_h,l_2]$. So, if the main term in Lemma \ref{lem-final} has been verified on shorter branching intervals, then
$$\M^{[l_1,l_2]}(\mu^{\D_h})\lessapprox 2^{\,2F_\mu(h)-F_\mu(j_h)-\min_{[h,j_h]}F_\mu}\cdot2^{2F_\mu(j_h)-\min_{[j_h+\delta k,l_2]}F_\mu-\min_{[j_h,l_2-\delta k]}F_\mu}.$$
If $j_h\geq h+\delta k$, due to $h\in [j_h-d_h, j_h]$, we have
$$F_\mu(j_h)=\min_{[h+\delta k, j_h]}F_\mu.$$
Since $F_\mu$ is Lipschitz,
$$F_\mu(j_h)-\min_{[h,j_h]}F_\mu=O(\delta) Tk,$$
that still holds when $j_h\leq h+\delta k$. Therefore, in either case, the exponent equals
$$2F_\mu(h)-\min_{[j_h+\delta k,l_2]}F_\mu-\min_{[j_h,l_2-\delta k]}F_\mu+O(\delta) Tk\leq 2F_\mu(h)-\min_{[h+\delta k,l_2]}F_\mu-\min_{[h,l_2-\delta k]}F_\mu+O(\delta) Tk.$$
Although it costs an additional $2^{O(\delta) Tk}$ and how much the branching interval shortens is not guaranteed, Case (M1) does not occur very often. If it occurs again at $j_{h'}^{[l_1', l_2']}>j_h^{[l_1, l_2]}$ with $l_2'-h'>10\delta k$, then $[l_1', l_2']\subset [l_1, l_2]$ and \eqref{jh-jjh} holds. If $j_{h'}-j_h< (l_2-h)/4$, then
$$(l_2-h)/4>d_{h'}-\delta k\geq (l_2'-h')/2-\delta k\geq (l_2'-h')/3.$$
This means the branching interval needs to be shortened for at least one quarter for another occurrence of Case (M1), thus $O(\log 1/\delta)$ occurrences in total. As it is the only case with the additional cost $2^{O(\delta) Tk}$, the accumulated error is $\lessapprox 2^{O(\delta\log1/\delta) Tk}$.

{\bf Case (M2).} $d_h<(l_2-h)/2$.

In this case, we apply Lemma \ref{lem-from-M-to-E} for
$$\M^{[l_1,l_2]}(\mu^{\D_h})\lessapprox 2^{-2(F_\mu(l_2-2d_h)-F_\mu(h))}\M^{[l_1,l_2]}(\mu^{\D_{l_2-2d_h}})+\max_{h<l\leq l_2-2d_h}2^{-2(F_\mu(l)-F_\mu(h))}\E_{l_2-l}(\mu^{\D_l}).$$
The reduction on $\E$ terms is already done. For $\M^{[l_1,l_2]}(\mu^{\D_{l_2-2d_h}})$, there are two cases.

{\bf Case (M2.1).} If $j_h> l_2-2d_h$, it goes to Case (M1).

{\bf Case (M2.2).} If $j_h\leq l_2-2d_h$, repeat the reduction on $\M$. It ends up with  either Case (M1) or a repetition of Case (M2.2). For the later one, by \eqref{jh-jjh}, the branching interval is shortened from $[h,l_2]$ by at least 
$$j_{l_2-2d_h}-j_h\geq d_{l_2-2d_h}-\delta k>4\delta k.$$
If the main term in Lemma \ref{lem-final} has been verified on shorter branching intervals, then
$$\begin{aligned}
\M^{[l_1,l_2]}(\mu^{\D_h})&\lessapprox 2^{-2(F_\mu(l_2-2d_{l_2-2d_h})-F_\mu(h))}\M^{[l_1,l_2]}(\mu^{\D_{l_2-2d_{l_2-2d_h}}})\\&\lessapprox 2^{2F_\mu(h)-\min_{[l_2-2d_{l_2-2d_h}+\delta k,l_2]}F_\mu-\min_{[l_2-2d_{l_2-2d_h},l_2-\delta k]}F_\mu}	
\\&\lessapprox 2^{2F_\mu(h)-\min_{[h+\delta k,l_2]}F_\mu-\min_{[h,l_2-\delta k]}F_\mu},
\end{aligned}
$$
as desired.

{\bf The accumulated error.} The reduction terminates when the length of the remaining branching interval is $<10\delta k$. There are two types of errors: one is $2^{O(\delta) Tk}$ from Case (M1) with $O(\log 1/\delta)$ occurrences; the other is $2^{O(\epsilon) Tk}$ with $O(1/\delta)$ occurrences because each time the branching interval is shortened for at least $\delta k$. In total, the accumulated error is
$$\lesssim 2^{O(\epsilon/\delta)Tk+O(\delta\log(1/\delta))Tk}$$
that completes the proof of Lemma \ref{lem-final}.
\end{proof}

\section{Proof of Proposition \ref{prop-good-bad}}\label{sec-good-bad}

In this section, we discuss the proof of Proposition \ref{prop-good-bad}. The planar version follows easily from arguments in \cite{Liu26+}\cite{GLW26+}, that relies on Orponen's $L^p$ radial projection estimates and classical $L^2$ orthogonal projection estimates. In general dimensions, we need the following radial projection theorem due to Ren (Theorem 1.13 in \cite{Ren23+}). Although it is not the same as what stated there, it follows directly from the argument.
\begin{prop}\label{prop-Ren}
Let $E_1, E_2\subset\R^d$ be compact separated sets, $d\ge2$, and $\mu_j$ be Frostman measures supported on $E_j$ of exponent $s>d/2$. If $d$ is odd, assume in addition that
every $\frac{d+1}{2}$-dimensional plane has measure zero.
Then, for every sufficiently small fixed $\eta>0$ and all $k\geq 1$, there are Borel sets
$\Gamma_k\subset E_1\times E_2$, with the following properties.
\begin{itemize}
 \item $\Gamma_k$ is decreasing, $
 (\mu_1\times\mu_2)(\Gamma_\infty)>0$ with $\Gamma_\infty:=\bigcap_k\Gamma_k$, $\mu_1\times\mu_2(\Gamma_k\backslash \Gamma_\infty)\lesssim 2^{-c_0\eta Tk}$.
 \item Each $\Gamma_k$ is a disjoint union of $O(2^{2d\eta Tk})$
products of Borel sets $\Gamma_k=\bigsqcup_iE_1^{k,i}\times E_2^{k,i}$.
 \item For every $\epsilon>0$, there are $c_\epsilon>0$, Borel sets
 $B_{k,\epsilon}\subset\Gamma_k$, such that
$$   (\mu_1\times\mu_2)(B_{k,\varepsilon})\lesssim 2^{-c_\epsilon Tk},$$
and $G_{k,\epsilon}:=\Gamma_k\setminus B_{k,\epsilon}$ satisfies that, for all $(x_1, x_2)\in E_1\times E_2$ and all
 $2^{-Tk}\le r\le1$,
$$\mu_1(T_{r, x_1, x_2}\cap\{y:(y,x_2)\in G_{k,\epsilon}\}), \,\mu_2(T_{r, x_1, x_2}\cap\{y:(x_1,y)\in G_{k,\epsilon}\})\leq  2^{\epsilon Tk}r^{\frac{d}{2}},$$
where $T_{r, x_1, x_2}$ denotes the $r\times \cdots \times r\times 1$ tube whose central line contains $x_1, x_2$.
\end{itemize}
Moreover, in even dimensions, one may take $\Gamma_k=E_1\times E_2$ for all $k\geq 1$.
\end{prop}

\begin{proof}[Proof of Proposition \ref{prop-good-bad}]
We only consider $(x_1, x_2)\in G_{k, \epsilon/2}$. Fix $x_1$ and denote
$$G_{k, \epsilon/2}^{x_1}:=\{x_2: (x_1, x_2)\in G_{k, \epsilon/2}\}.$$
We first consider tubes. For $0\leq l<m\leq k$, let $T(x_1,\omega)\subset Q_l(x_1)$ denote $2^{-Tm}\times\cdots\times 2^{-Tm}\times 2^{-Tl}$ tubes containing $x_1$ parallel to $\omega$. By Markov's inequality,
$$\begin{aligned}& \mu_2\{x_2\in G_{k, \epsilon/2}^{x_1}: T(x_1,\omega_{x_1, x_2})\text{ is bad}\}\\\lesssim &2^{-\epsilon Tk}\cdot 2^{\frac{d}{2}T(m-l)+\min_{[l,m]} (F_{\mu}-F_{\mu}(l))}\int_{G_{k, \epsilon/2}^{x_1}} (\mu_1)_{Q_l(x_1)}(T(x_1,\omega_{x_1, x_2}))d\mu_2(x_2).\end{aligned}$$
Decompose $S^{d-1}$ into $2^{-T(m-l)}$-caps. Then
$$\begin{aligned}&\int_{G_{k, \epsilon/2}^{x_1}} (\mu_1)_{Q_l(x_1)}(T(x_1,\omega_{x_1, x_2}))d\mu_2(x_2)\\\approx &\sum_{\omega: 2^{-T(m-l)}\text{-caps}}(\mu_1)_{Q_l(x_1)}(T(x_1,\omega_{x_1, x_2}))\cdot\mu_2\{x_2\in G_{k, \epsilon/2}^{x_1}:\omega_{x_1, x_2}\in\omega\},\end{aligned}$$
which is, by Proposition \ref{prop-Ren} and counting the overlapping between rotated tubes,
$$\lesssim 2^{\epsilon Tk/2}\cdot2^{-\frac{d}{2}T(m-l)}\sum_{j=l}^m2^{\frac{d}{2}T(j-l)} (\mu_1)_{Q_l}(Q_j)\lesssim  (m-l)2^{\epsilon Tk/2-\frac{d}{2}T(m-l)-\min_{[l,m]} (F_{\mu}-F_{\mu}(l))}.$$

Overall,
$$\mu_2\{x_2\in G_{k, \epsilon/2}^{x_1}: T(x_1,\omega_{x_1, x_2})\text{ is bad}\}\lesssim (m-l)2^{-\epsilon Tk/2},$$
so the good tube estimate in Proposition \ref{prop-good-bad} follows from removing all bad pairs $(x_1, x_2)$ over $0\leq l<m\leq k$.

Now we turn to slabs. For $0\leq l<m\leq k$, let $S(x,\omega)\subset Q_l(x_1)$ denote $2^{-Tl}\times\cdots\times 2^{-Tl}\times 2^{-Tm}$ slabs containing $x$ perpendicular to $\omega$. By Markov's inequality,
\begin{equation}\label{slab-Markov}
	\begin{aligned}& \mu_2\{x_2\in G_{k, \epsilon/2}^{x_1}: S(x,\omega_{x_1, x_2})\text{ is bad on average in }Q_l(x_1)\}\\\lesssim &2^{-\epsilon Tk}\cdot 2^{T(m-l)+\min_{[l,m]} (F_{\mu}-F_{\mu}(l))}\int_{G_{k, \epsilon/2}^{x_1}} \int (\mu_1)_{Q_l(x_1)}(S(x,\omega_{x_1, x_2}))d(\mu_1)_{Q_l(x_1)}(x)d\mu_2(x_2).\end{aligned}
\end{equation}
If we let $\pi^{x_1}(x_2):=\omega_{x_1,x_2}$ denote the radial projection and $\nu:=\pi^{x_1}_*(\mu_2|_{G_{k, \epsilon/2}^{x_1}})$ be the pushforward measure on $S^{d-1}$, then Proposition \ref{prop-Ren} implies that $\nu$ is a Frostman measure of exponent $d/2$ up to a factor $2^{\epsilon Tk/2}$. Therefore, the integral in \eqref{slab-Markov} is equivalent to
$$\iint (\mu_1)^{Q_l(x_1)}(2^{Tl}S(x,\omega))d(\mu_1)^{Q_l(x_1)}(x)d\nu(\omega)\lesssim 2^{-T(m-l)}\int\| \pi_\omega ((\mu_1)^{Q_l}*\psi_{m-l})\|_{L^2(\R \omega)}^2d\nu(\omega).$$
where $\pi_\omega$ denotes the orthogonal projection to $\R\omega$ and $\psi_{m-l}(\cdot)=2^{dT(m-l)}\psi(2^{T(m-l)}\cdot)$ is a nonnegative mollifier with bounded Fourier support. Then, by a classical $L^2$ argument on orthogonal projections (see, e.g., the proof of Theorem 5.6 in \cite{Mat15}), it is
$$\begin{aligned}
&\lesssim 2^{-T(m-l)}\cdot 2^{\epsilon Tk/2}\int|\widehat{(\mu_1)^{Q_l}}(\xi)\hat{\psi}(2^{-T(m-l)}\xi)|^2(1+|\xi|)^{-d/2}d\xi\\&\lesssim 2^{-T(m-l)}\cdot 2^{\epsilon Tk/2}\sum_{j=0}^{m-l} 2^{-\frac{d}{2}Tj}\int|\widehat{(\mu_1)^{Q_l}}(\xi)|^2\psi(2^{-Tj}\xi)d\xi\\&\leq 2^{-T(m-l)}\cdot 2^{\epsilon Tk/2}\sum_{j=0}^{m-l} 2^{\frac{d}{2}Tj}\iint|\hat{\psi}(2^{Tj}(x-y)|\,d(\mu_1)^{Q_l}(x)\,d(\mu_1)^{Q_l}(y)\\&\lesssim (m-l)2^{-T(m-l)}\cdot 2^{\epsilon Tk/2}\cdot 2^{-\min_{[l,m]} (F_{\mu}-F_{\mu}(l))}.
\end{aligned}$$
Overall,
$$\mu_2\{x_2\in G_{k, \epsilon/2}^{x_1}: S(x,\omega_{x_1, x_2})\text{ is bad on average in }Q_l(x_1)\}\lesssim (m-l)2^{-\epsilon Tk/2},$$
so the good slab estimate in Proposition \ref{prop-good-bad} follows from removing all bad pairs $(x_1, x_2)$ over $0\leq l<m\leq k$.
\end{proof}

\bibliographystyle{plain}
\bibliography{mybibtex.bib}

\end{document}